\documentclass[11pt]{article}
\usepackage{latexsym,amssymb}
\usepackage{amsmath}
\usepackage[mathscr]{eucal}
\usepackage{color}
\usepackage[abbrev]{amsrefs}

\newtheorem{Theorem}{\bf Theorem}[section]
\newtheorem{Lemma}{\bf Lemma}[section]
\newtheorem{Proposition}{\bf Proposition}[section]
\newtheorem{Corollary}{\bf Corollary}[section]
\newtheorem{Remark}{\bf Remark}[section]
\newtheorem{Example}{\bf Example}[section]
\newtheorem{Definition}{\bf Definition}[section]

\newenvironment{theorem}{\begin{Theorem}$\!\!\!$}{\end{Theorem}}
\newenvironment{lemma}{\begin{Lemma}$\!\!\!$}{\end{Lemma}}
\newenvironment{proposition}{\begin{Proposition}$\!\!\!$}{\end{Proposition}}
\newenvironment{corollary}{\begin{Corollary}$\!\!\!$}{\end{Corollary}}
\newenvironment{remark}{\begin{Remark}$\!\!\!$}{\end{Remark}}

\newenvironment{definition}{\begin{Definition}$\!\!\!$}{\end{Definition}}

\numberwithin{equation}{section}

\numberwithin{equation}{section}

\newcommand{\dee}{{\rm{d}}}

\def\Xint#1{\mathchoice
{\XXint\displaystyle\textstyle{#1}}%
{\XXint\textstyle\scriptstyle{#1}}%
{\XXint\scriptstyle\scriptscriptstyle{#1}}%
{\XXint\scriptscriptstyle\scriptscriptstyle{#1}}%
\!\int}
\def\XXint#1#2#3{{\setbox0=\hbox{$#1{#2#3}{\int}$}
\vcenter{\hbox{$#2#3$}}\kern-.5\wd0}}

\def\dashint{\Xint-}

\allowdisplaybreaks

\begin{document}

\title{Initial traces and solvability of the fast diffusion equation\\ with power-type nonlinearity}
\author{
\qquad\\
Kazuhiro Ishige and Nobuhito Miyake
}
\date{}
\maketitle
\begin{abstract}
We investigate the qualitative behavior of the initial traces of nonnegative solutions to the fast diffusion equation with power-type nonlinearity.
Necessary conditions for the existence of solutions to the corresponding Cauchy problem are identified.
Moreover, sharp sufficient conditions for the existence of solutions are established by employing uniform local Morrey spaces and their generalizations.
\end{abstract}
\vspace{25pt}
\noindent Addresses:

\smallskip
\noindent 
K.~I.:  Graduate School of Mathematical Sciences, The University of Tokyo,\\
\qquad\,\,\, 3-8-1 Komaba, Meguro-ku, Tokyo 153-8914, Japan. \\
\noindent 
E-mail: {\tt ishige@ms.u-tokyo.ac.jp}\\

\smallskip
\noindent 
N.~M.:  Faculty of Mathematics, Kyushu University, Fukuoka, 819-0395, Japan.\\
\noindent 
E-mail: {\tt miyake@math.kyushu-u.ac.jp}\\
\date{}
\maketitle

\vspace{25pt}
\vspace{10pt}
\vspace{3pt}

\newpage
%%%%%%%%%%%%%%%%%%%%%%%
%%%%%%%%%%%%%%%%%%%%%%
\section{Introduction}
%%%%%%%%%%%%%%%%%%%%%%
%%%%%%%%%%%%%%%%%%%%%%
In this paper, we study the qualitative behavior of the initial traces of nonnegative solutions to 
the fast diffusion equation with a power-type nonlinearity:
\begin{equation}
\tag{E}
\label{eq:E}
\partial_t u=\Delta u^m+u^p\quad\quad\mbox{in}\quad{\mathbb R}^N\times(0,T),
\end{equation}
where $N\ge 1$, $T\in(0,\infty]$, $m\in(0,1)$, and $p>1$. 
We establish necessary conditions for the existence of nonnegative solutions to the Cauchy problem
\begin{equation}
\tag{P}
\label{eq:P}
\left\{
\begin{array}{ll}
\partial_t u=\Delta u^m+u^p & \quad\mbox{in}\quad{\mathbb R}^N\times(0,T),\vspace{3pt}\\
u(\cdot,0)=\mu & \quad\mbox{in}\quad{\mathbb R}^N,
\end{array}
\right.
\end{equation}
where $\mu$ is a nonnegative Radon measure on ${\mathbb R}^N$. 
Furthermore, we establish sharp sufficient conditions for the existence of nonnegative solutions 
to problem~\eqref{eq:P}, and identify optimal singular profiles of the initial data 
for which problem~\eqref{eq:P} is solvable. 

Throughout this paper, 
$\mathcal{M}$ denotes the set of nonnegative Radon measures on~$\mathbb{R}^N$,
and $\mathcal{L}$ denotes the set of nonnegative locally integrable functions on~$\mathbb{R}^N$.
We often identify $\dee\mu = \mu(x)\,\dee x$ in $\mathcal{M}$ for $\mu \in \mathcal{L}$.
Let ${\mathcal L}^N$ be the $N$-dimensional Lebesgue measure.
For any $f\in\mathcal{L}$, $z\in{\mathbb R}^N$, and $r>0$, 
we define the average of $f$ over the ball $B(z, r)$ by
$$
\dashint_{B(z,r)} f\,\dee x:=\frac{1}{{\mathcal L}^N(B(z,r))}\int_{B(z,r)} f\,\dee x,
$$
where $B(z,r):=\{x\in{\mathbb R}^N\,:\,|x-z|<r\}$. 
We also introduce the constants
$$
p_m:=m+\frac{2}{N},\quad \theta:=\frac{p-m}{2(p-1)},\quad
\theta':=\frac{1}{\theta}=\frac{2(p-1)}{p-m},\quad \kappa:=N(m-1)+2.
$$
Note that $p_m>1$ if and only if $\kappa>0$. 
Furthermore, we define
\begin{equation}
\label{eq:1.1}
\Psi_\alpha(\xi):=\xi[\log (e+\xi)]^\alpha,
\qquad
\eta(\xi):=\xi^N\biggr[\log\biggr(e+\frac{1}{\xi}\biggr)\biggr]^{\frac{N}{2}},
\end{equation}
for $\xi \in [0, \infty)$, where $\alpha>0$ is a fixed parameter. 

The study of initial traces of solutions to parabolic equations is useful for deriving necessary conditions for the existence of solutions 
to the corresponding Cauchy problem. 
This subject has attracted considerable attention from many mathematicians; 
see, for instance, \cites{A, HW, Hui, W1, W2} for linear parabolic equations, 
\cites{AC, HP} for porous medium equations, 
\cites{DH, DH02} for parabolic $p$\,-Laplace equations, 
\cites{I, IJK, ZX} for doubly nonlinear parabolic equations,
\cite{BSV} for fractional diffusion equations,
\cite{AIS} for Finsler heat equations,
\cites{ADi, BP, FHIL, FHIL02, Hisa, HI18, HI24, HIT02, FI01, IKO, IMS, TY} for parabolic equations with source (positive) nonlinearity, 
\cites{BCV, BD, MV, MV02, MV03} for parabolic equations with absorption (negative) nonlinearity.
Among these, the qualitative behavior of the initial traces of nonnegative solutions to equation~\eqref{eq:E} with $1 \le m < p$
has been investigated in \cites{ADi, IMS}. 
It has been shown that the initial traces of such solutions exhibit the following properties.
\begin{itemize}
  \item[(A)]
  Let $1\le m<p$.
  Assume that equation~\eqref{eq:E} admits a solution~$u$ in ${\mathbb R}^N\times(0,T)$ for some $T\in(0,\infty)$. 
  Then there exists a unique $\nu\in{\mathcal M}$ such that 
  $$
  \underset{t\to+0}{\mbox{{\rm ess lim}}}\int_{{\mathbb R}^N} u(t)\psi\,\dee x=\int_{{\mathbb R}^N}\psi\,\dee \nu,\quad \psi\in C_c({\mathbb R}^N).
  $$
  Furthermore, there exists $C=C(N,m,p)>0$ such that  
  \begin{equation*}
  \left\{
  \begin{array}{ll}
  \displaystyle{\sup_{z\in{\mathbb R}^N}\nu(B(z,T^\theta))\le CT^{\theta\left(N-\frac{2}{p-m}\right)}} & \mbox{if $p<p_m$},\\
   \displaystyle{\sup_{z\in{\mathbb R}^N}\nu(B(z,\sigma))\le C\left[\log\left(e+\frac{T^{\theta}}{\sigma}\right)\right]^{-\frac{N}{2}}} 
    & \mbox{for $\sigma\in(0,T^\theta)$ if $p=p_m$},\vspace{5pt}\\
  \displaystyle{\sup_{z\in{\mathbb R}^N}\nu(B(z,\sigma))\le C\sigma^{N-\frac{2}{p-m}}}
    & \mbox{for $\sigma\in(0,T^\theta)$ if $p>p_m$}.
  \end{array}
  \right.
  \end{equation*}
  \item[(B)] 
  Let $u$ be a solution to problem~\eqref{eq:P} with $1\le m<p$. 
  Then $u$ is a solution to equation~\eqref{eq:E} and 
  the initial trace of $u$ coincides with the initial data of $u$ in ${\mathcal M}$.
\end{itemize} 
In view of assertion~(B), assertion~(A) provides necessary conditions for the existence of solutions to problem~\eqref{eq:P}.
On the other hand, in \cite{IMS}, the present authors, together with Sato,
established the following sufficient conditions for the existence of solutions to problem~\eqref{eq:P}
(see also \cite{ADi}*{Theorem~3.1} for assertion~(C)-(1)).
\begin{itemize}
\item[(C)] 
 Let $1\le m<p$. 
  \begin{itemize}
  \item[(1)]
  Let $p<p_m$. 
  Then problem~\eqref{eq:P} admits a local-in-time solution if and only if 
  $\displaystyle{\sup_{z\in{\mathbb R}^N}\mu(B(z,1))<\infty}$. 
  \item[(2)] 
  Let $p=p_m$. Then, for any $\alpha>0$, there exists $\delta_1=\delta_1(N,m,\alpha)>0$ such that 
  if $\mu\in\mathcal{L}$ satisfies
  $$
  \sup_{z\in{\mathbb R}^N}\,\sup_{\sigma\in(0,T^{\theta})}\,
  \eta\left(\frac{\sigma}{T^{\theta}}\right)
  \Psi_\alpha^{-1}\left(\,\dashint_{B(z,\sigma)}\Psi_\alpha\left(T^{\frac{1}{p-1}}\mu\right)\,\dee x\right)
  \le\delta_1
 $$
  for some $T\in(0,\infty)$,
  then problem~\eqref{eq:P} admits a solution in ${\mathbb R}^N\times(0,T)$. 
  Here, $\Psi_\alpha$ and $\eta$ are defined as in \eqref{eq:1.1}.  
  \item[(3)] 
  Let $p>p_m$. 
  Then, for any $\beta>1$, there exists $\delta_2=\delta_2(N,m, p,\beta)>0$ such that 
  if $\mu\in\mathcal{L}$ satisfies
  $$
  \sup_{z\in{\mathbb R}^N}\sup_{\sigma\in(0,T^\theta)}\sigma^{\frac{2}{p-m}}\left(\,\dashint_{B(z,\sigma)}\mu^\beta\,\dee x\right)^{\frac{1}{\beta}}\le \delta_2
  $$
  for some $T\in(0,\infty]$, then problem~\eqref{eq:P} admits a solution in ${\mathbb R}^N\times(0,T)$.
  \end{itemize}
\end{itemize}
Combining the results in assertions (A), (B), and~(C), 
we have the following assertion.
\begin{itemize}
\item[(D)]
Let $m\ge 1$ and $p\ge p_m$. 
For any $c>0$, set 
$$
\mu_c(x):=
\left\{
\begin{array}{ll}
c|x|^{-N}\displaystyle{\biggr[\log\biggr(e+\frac{1}{|x|}\biggr)\biggr]^{-\frac{N}{2}-1}} & \mbox{if}\quad p=p_m,\vspace{7pt}\\
c|x|^{-\frac{2}{p-m}} & \mbox{if}\quad p>p_m,
\end{array}
\right.
$$
for almost all~$x\in{\mathbb R}^N$. 
\begin{itemize}
  \item[{\rm (1)}] 
  Problem~\eqref{eq:P} admits a local-in-time solution for sufficiently small $c>0$; 
  \item[{\rm (2)}] 
  Problem~\eqref{eq:P} admits no local-in-time solutions for sufficiently large $c>0$. 
\end{itemize}
Furthermore, if $p>p_m$ and $c>0$ is sufficiently small, then 
problem~\eqref{eq:P} admits a global-in-time solution.
\end{itemize}
The results in (D) show that the strength of the singularity at the origin of the functions $\mu_c$
serves as the critical threshold for the local solvability of problem~\eqref{eq:P}.
They also identify optimal singular profiles of the initial data for which problem~\eqref{eq:P} is solvable when $m\ge 1$, 
thereby demonstrating the sharpness of the sufficient conditions given in assertion~(C).

We also recall sufficient conditions for the existence of solutions to problem~\eqref{eq:P} 
with measures as initial data.
\begin{itemize}
\item[(E)]
Let $m\ge1$ and $p>p_m$.
Then problem~\eqref{eq:P} admits a local-in-time solution if $\mu\in\mathcal{M}$ satisfies
$$
\sup_{z\in\mathbb{R}^N}\sup_{\sigma\in(0,1)}\sigma^{\gamma-N}\mu(B(z, \sigma))<\infty
$$
for some $\gamma\in(0,2/(p-m))$ (see \cite{An}*{Remark~1.2};  see also \cite{Ta} for further results in the case of $m=1$.) 
\end{itemize}
However, assertion~(E) does not apply to the function $\mu_c$ in assertion~(D).

In contrast, in the fast diffusion case $m\in(0,1)$, as far as we know, 
no results are available on the qualitative behavior of the initial traces of solutions to problem~\eqref{eq:P}, 
nor on the identification of optimal singular profiles of the initial data for which problem~\eqref{eq:P} is solvable.

The aim of this paper is to establish analogous results to assertions~(A), (B), (C), and (D) 
in the fast diffusion case $m\in(0,1)$.
These results extend previous works such as~\cites{IMS, ADi} to the case $m\in(0,1)$.
Our approach follows the arguments in~\cites{IMS, ADi}, 
but requires a more refined analysis to overcome the additional difficulties arising in the case $m\in(0,1)$.
\vspace{5pt}

We introduce the definitions of solutions to equation~\eqref{eq:E} and problem~\eqref{eq:P}, respectively.
\begin{definition}
\label{Definition:1.1}
Let $N\ge 1$, $m\in(0,1)$, $p>1$, $T\in(0,\infty]$, and $u\in L^p_{\rm loc}(\mathbb{R}^N\times[0, T))$ be nonnegative in ${\mathbb R}^N\times(0,T)$. 
\begin{itemize}
  \item[{\rm (1)}] 
  We say that $u$ is a solution to 
  equation~\eqref{eq:E} in ${\mathbb R}^N\times(0,T)$ if $u$ satisfies 
  $$
  \int^T_\tau\int_{\mathbb{R}^N}(-u\partial_t \phi-u^m\Delta\phi-u^p\phi)\,\dee x\,\dee t=\int_{\mathbb{R}^N}u(x,\tau)\phi(x,\tau)\,\dee x
  $$
  for $\phi\in C^{2;1}_c({\mathbb R}^N\times[0,T))$ and almost all $\tau\in(0,T)$.
  \item[{\rm (2)}] 
  For any $\mu\in{\mathcal M}$,  
  we say that $u$ is a solution to problem~\eqref{eq:P} in ${\mathbb R}^N\times(0,T)$ if $u$ satisfies 
  $$
  \int^T_0 \int_{\mathbb{R}^N}(-u\partial_t \phi-u^m\Delta\phi-u^p\phi)\,\dee x\,\dee t=\int_{\mathbb{R}^N}\phi(x,0)\,\dee\mu(x)
  $$
  for $\phi\in C^{2;1}_c({\mathbb R}^N\times[0,T))$. 
\end{itemize}
\end{definition}

We introduce some notation. 
For any positive functions $f$ and $g$ in a set $X$, 
we say that $f(x)\preceq g(x)$ for $x\in X$ or equivalently that $g(x)\succeq f(x)$ for $x\in X$ if 
there exists $C>0$ such that 
$$
f(x)\le Cg(x)\quad\mbox{for $x\in X$}.
$$
If $f(x)\preceq g(x)$ and $g(x)\preceq f(x)$ for $x\in X$, we say that $f(x)\asymp g(x)$ for $x\in X$. 
Throughout this paper, we use $C$ to denote generic positive constants, which may vary from line to line.

Let $\Phi$ be a nonnegative, convex, and strictly increasing function in $[0,\infty)$ such that $\Phi(0)=0$. 
Let $\rho$ be a nonnegative, non-decreasing, and continuous function in $[0,\infty)$. 
Then, for any $f\in\mathcal{L}$ and $R\in(0,\infty]$, 
define
$$
|||f|||_{\rho,\Phi;R}:=\sup_{z\in{\mathbb R}^N}\sup_{\sigma\in(0,R)}
\left\{\rho(\sigma)\Phi^{-1}\left(\,\dashint_{B(z,\sigma)}\Phi(f)\,\dee x\right)\right\}.
$$
For any $q\in[1,\infty)$ and $\alpha\in[1,\infty)$, 
if $\Phi(\xi)=\xi^\alpha$ and $\rho(\xi)=\xi^{N/q}$ for $\xi\in[0,\infty)$, 
we write 
$$
|||f|||_{q,\alpha;R}:=|||f|||_{\rho,\Phi;R}
=\sup_{z\in{\mathbb R}^N}\sup_{\sigma\in(0,R)}
\left\{\sigma^{\frac{N}{q}}\left(\,\dashint_{B(z,\sigma)}|f|^\alpha\,\dee x\right)^{\frac{1}{\alpha}}\right\}
$$
for simplicity. 
Note that 
$$
|||f|||_{q,\alpha,R}<\infty\quad\mbox{if and only if}\quad |||f|||_{q,\alpha,1}<\infty
$$
for $q\in[1,\infty)$, $\alpha\in[1,\infty)$, and $R\in(0, \infty)$. 
This equivalence follows from the following fact:
for any fixed $k>1$, there exists $C>0$ such that 
\begin{equation}
\label{eq:1.2}
 \sup_{z\in{\mathbb R}^N}\nu({B(z,kR))\le C\sup_{z\in{\mathbb R}^N}\nu(B(z,R)})
\end{equation}
for $\nu\in{\mathcal M}$ and $R>0$.
\vspace{5pt}

We are now in a position to state the main results of this paper. 
The first main result concerns the qualitative behavior of the initial traces of solutions to equation~\eqref{eq:E}, 
and provides necessary conditions for the existence of solutions to problem~\eqref{eq:P}.%
\begin{theorem}
\label{Theorem:1.1}
Let $N\ge 1$, $m\in(0,1)$, and $p>1$. 
\begin{itemize}
  \item[{\rm (1)}] 
  Let $u$ be a solution to equation~\eqref{eq:E} in ${\mathbb R}^N\times(0,T)$, where $T\in(0,\infty)$. 
  Then there exists a unique $\nu\in {\mathcal M}$ such that 
  $$
  \underset{t\to+0}{\mbox{{\rm ess lim}}}\int_{{\mathbb R}^N}u(t)\psi\,\dee x=\int_{{\mathbb R}^N}\psi\,\dee\nu(x),
  \quad \psi\in C_c({\mathbb R}^N).
  $$
  Furthermore, there exists $C=C(N,m,p)>0$ such that
  \begin{equation*}
  \left\{
  \begin{array}{ll}
  \displaystyle{\sup_{z\in{\mathbb R}^N}\nu(B(z,T^\theta))\le CT^{\theta\left(N-\frac{2}{p-m}\right)}} & \mbox{if $1<p<p_m$},\\
   \displaystyle{\sup_{z\in{\mathbb R}^N}\nu(B(z,\sigma))\le C\left[\log\left(e+\frac{T^{\theta}}{\sigma}\right)\right]^{-\frac{N}{2}}} 
    & \mbox{for $\sigma\in(0,T^\theta)$ if $p=p_m$},\vspace{5pt}\\
  \displaystyle{\sup_{z\in{\mathbb R}^N}\nu(B(z,\sigma))\le C\sigma^{N-\frac{2}{p-m}}}
    & \mbox{for $\sigma\in(0,T^\theta)$ if $p>p_m$}.
  \end{array}
  \right.
  \end{equation*}
  \item[{\rm (2)}] 
  Let $u$ be a solution to problem~\eqref{eq:P}. 
  Then $u$ is a solution to equation~\eqref{eq:E} and 
  the initial trace of $u$ coincides with the initial data of $u$ in ${\mathcal M}$.
\end{itemize}
\end{theorem}
\begin{remark}
\label{Remark:1.1}
{\rm (1)} Let $u$ be a solution to equation~\eqref{eq:E} in ${\mathbb R}^N\times(0,T)$, where $T\in(0,\infty]$. 
For any $\lambda>0$, define 
$$
u_\lambda(x,t):=\lambda^{\frac{2}{p-m}}u(\lambda x,\lambda^{\theta'} t),\quad (x,t)\in{\mathbb R}^N\times(0,T_\lambda),
$$
where $T_\lambda:=\lambda^{-\theta'}T$. Then $u_\lambda$ is a solution to equation~\eqref{eq:E} in ${\mathbb R}^N\times(0,T_\lambda)$.
\vspace{3pt}
\newline
{\rm (2)} 
If $1<p\le p_m$, then problem~\eqref{eq:P} admits no nontrivial global-in-time solutions 
{\rm ({\it see} \cites{MM,Q})}. 
This result follows immediately from Theorem~{\rm\ref{Theorem:1.1}}.
Indeed, assume that problem~\eqref{eq:P} admits a nontrivial global-in-time solution $u$.
Then, for almost all~$\tau>0$, define 
$$
u_\tau(x,t):=u(x,t+\tau)\quad\mbox{for almost all~$(x,t)\in{\mathbb R}^N\times(0,\infty)$}.
$$ 
Then $u_\tau$ is a global-in-time solution to problem~\eqref{eq:P} with initial data $\mu=u(\tau)$.
Applying Theorem~{\rm\ref{Theorem:1.1}} yields:  
\begin{align*}
 & \sup_{z\in{\mathbb R}^N}\int_{B(z,T^{\theta})} u(\tau)\,\dee x\le CT^{N\theta-\frac{1}{p-1}}\to 0\quad\mbox{if}\quad p<p_m,\\
 & \sup_{z\in{\mathbb R}^N}\int_{B(z,T^{\frac{\theta}{2}})} u(\tau)\,\dee x
\le C\left[\log\left(e+T^{\frac{\theta}{2}}\right)\right]^{-\frac{N}{2}}\to 0\quad\mbox{if}\quad p=p_m,
\end{align*}
as $T\to\infty$. This implies that $u(x,\tau)=0$ for almost all~$(x,\tau)\in{\mathbb R}^N\times(0,\infty)$, 
contradicting the assumption that $u$ is nontrivial.
Hence, problem~\eqref{eq:P} admits no nontrivial global-in-time solutions when $1<p\le p_m$.
\end{remark}
We next state our main results on sufficient conditions 
for the existence of solutions to problem~\eqref{eq:P}. 
Theorem~\ref{Theorem:1.2} shows that 
the necessary condition given in Theorem~\ref{Theorem:1.1} is also sufficient in the subcritical case $1<p<p_m$. 
This theorem corresponds to assertion~(C)-(1) in the case $m\ge 1$. 
\begin{theorem}
\label{Theorem:1.2}
Let $N\ge 1$, $m\in(0,1)$, and $1<p<p_m$. 
\begin{itemize}
  \item[{\rm (1)}]  
  There exists $\delta>0$ such that if $\mu\in\mathcal{M}$ satisfies
  \begin{equation}
  \label{eq:1.3}
  \sup_{z\in\mathbb{R}^N}\mu(B(z,T^\theta))\le\delta T^{\theta\left(N-\frac{2}{p-m}\right)}
  \end{equation}
  for some $T\in(0,\infty)$, 
  then problem~\eqref{eq:P} admits a solution in ${\mathbb R}^N\times(0,T)$, with $u$ satisfying 
  $$
  T^{\frac{1}{p-1}}\sup_{t\in(0,T)} \sup_{z\in\mathbb{R}^N}\dashint_{B(z, T^\theta)}u(t)\,\dee x
  +
  T^{-\frac{N}{\kappa}+\frac{1}{p-1}}\sup_{t\in(0,T)}\left\{t^{\frac{N}{\kappa}}\|u(t)\|_{L^\infty({\mathbb R}^N)}\right\}\le C_*,
  $$
  where $C_*$ is a constant depending only on $N$, $m$, and $p$.
  \item[{\rm (2)}]
  Problem~\eqref{eq:P} admits a local-in-time solution with initial data $\mu\in\mathcal{M}$ if and only if 
  $\displaystyle{\sup_{z\in{\mathbb R}^N}\mu(B(z,1))<\infty}$.
\end{itemize}
\end{theorem}

In Theorems~\ref{Theorem:1.3} and \ref{Theorem:1.4}, 
we provide sufficient conditions for the existence of solutions to problem~\eqref{eq:P} 
in the critical case $p=p_m$ and the supercritical case $p>p_m$, respectively, 
which make it possible to identify the optimal singular profile.
These theorems correspond, respectively, to assertions~(C)-(2) and~(C)-(3) in the case $m\ge 1$.
\begin{theorem}
\label{Theorem:1.3}
Let $N\ge 1$, $m\in(0,1)$, $p=p_m>1$, and $\alpha>0$. 
Let $\Psi_\alpha$ and $\eta$ be defined as in \eqref{eq:1.1}. 
Then there exists $\delta>0$ such that
if $\mu\in\mathcal{L}$ satisfies
\begin{equation}
\label{eq:1.4}
\sup_{z\in{\mathbb R}^N}\,\sup_{\sigma\in(0,T^{\theta})}\,
\left\{
\eta\left(\frac{\sigma}{T^{\theta}}\right)
\Psi_\alpha^{-1}\left(\,\dashint_{B(z,\sigma)}\Psi_\alpha\left(T^{\frac{1}{p-1}}\mu\right)\,\dee x\right)
\right\}
\le\delta
\end{equation}
for some $T\in(0,\infty)$,
then problem~\eqref{eq:P} admits a solution $u$ in ${\mathbb R}^N\times(0,T)$, with $u$ satisfying
\begin{align*}
 & \sup_{t\in(0,T)}\sup_{z\in{\mathbb R}^N}\,\sup_{\sigma\in(0,T^{\theta})}\,
\left\{
\eta\left(\frac{\sigma}{T^{\theta}}\right)
\Psi_\alpha^{-1}\left(\,\dashint_{B(z,\sigma)}\Psi_\alpha\left(T^{\frac{1}{p-1}}u(t)\right)\,\dee x\right)
\right\}\\
 & \qquad\qquad\qquad+\sup_{t\in(0,T)}
 \left\{t^{\frac{1}{p-1}}\left[\log\left(e+\frac{T}{t}\right)\right]^{\frac{1}{p-1}}\|u(t)\|_{L^\infty({\mathbb R}^N)}\right\}\le C_*,
\end{align*}
where $C_*$ is a constant depending only on $N$, $m$, $p$, and $\alpha$.
\end{theorem}
\begin{theorem}
\label{Theorem:1.4}
Let $N\ge 1$, $m\in(0,1)$, $p>1$ with $p>p_m$. 
Let $\beta\in(1,N(p-m)/2)$ be such that $\kappa_\beta:=N(m-1)+2\beta>0$. 
Then there exists $\delta>0$ such that
if $\mu\in\mathcal{L}$ satisfies
\begin{equation}
\label{eq:1.5}
|||\mu|||_{\frac{N(p-m)}{2},\beta;\,T^{\theta}}\le\delta
\end{equation}
for some $T\in(0,\infty]$, 
then problem~\eqref{eq:P} admits a solution $u$ in ${\mathbb R}^N\times(0,T)$, with $u$ satisfying
\begin{equation}
\label{eq:1.6}
\sup_{t\in(0,T)}|||u(t)|||_{\frac{N(p-m)}{2},\beta;\,T^{\theta}}
+\sup_{t\in(0,T)}\left\{t^{\frac{1}{p-1}}\|u(t)\|_{L^\infty({\mathbb R}^N)}\right\}\le C_*,
\end{equation}
where $C_*$ is a constant depending only on $N$, $m$, $p$, and $\beta$.
\end{theorem}
We emphasize that Theorem~\ref{Theorem:1.4} provides a criterion for initial data to ensure the existence of global-in-time solutions to problem~\eqref{eq:P}.

\begin{Remark}
\label{Remark:1.2}
It was shown in \cite{Ta} that, 
in the case $m=1$ with $p=p_1$ {\rm({\it resp.~$p>p_1$})}, 
there exists $\mu\in \mathcal{L}$ satisfying \eqref{eq:1.4} with $\alpha=0$ {\rm({\it resp.~\eqref{eq:1.5} with $\beta=1$})} 
such that problem~\eqref{eq:P} admits no local-in-time solutions.
This suggests that 
Theorem~{\rm\ref{Theorem:1.3}} with $\alpha=0$ and Theorem~{\rm\ref{Theorem:1.4}} with $\beta=1$ do not hold.
\end{Remark}

Sufficient conditions given in Theorems~\ref{Theorem:1.3} and \ref{Theorem:1.4} are sharp. 
Indeed, similarly to assertion~(D), 
these results, together with Theorem~\ref{Theorem:1.1}, enable us to identify the optimal singular profiles 
of the initial data for which problem~\eqref{eq:P} with $p\ge p_m$ is solvable. 
\begin{corollary}
\label{Corollary:1.1}
Let $N\ge 1$, $m\in(0,1)$, and $p>1$ with $p\ge p_m$. 
For any $c>0$, set 
$$
\mu_c^m(x):=
\left\{
\begin{array}{ll}
c|x|^{-N}\displaystyle{\biggr[\log\biggr(e+\frac{1}{|x|}\biggr)\biggr]^{-\frac{N}{2}-1}} & \mbox{if}\quad \displaystyle{p=p_m},\vspace{7pt}\\
c|x|^{-\frac{2}{p-m}} & \mbox{if}\quad \displaystyle{p>p_m},\vspace{3pt}
\end{array}
\right.
$$
for almost all~$x\in{\mathbb R}^N$. 
\begin{itemize}
  \item[{\rm (1)}] 
  Problem~\eqref{eq:P} admits a local-in-time solution for sufficiently small $c>0$; 
  \item[{\rm (2)}] 
  Problem~\eqref{eq:P} admits no local-in-time solutions for sufficiently large $c>0$. 
\end{itemize}
Furthermore, if $p>p_m$ and $c>0$ is sufficiently small, then 
problem~\eqref{eq:P} admits a global-in-time solution.
\end{corollary}

The proofs of our theorems are based on the analytical framework developed in \cites{ADi, IMS}.
In particular, the proof of Theorem~\ref{Theorem:1.1} is a slight modification of the argument presented in \cite{IMS} (see also \cite{TY}).
In contrast, the proofs of Theorems~\ref{Theorem:1.3} and \ref{Theorem:1.4} require more delicate analysis. 
The proof of Theorem~\ref{Theorem:1.2} is an application of the 
methods developed in \cites{ADi} and the proofs of Theorems~\ref{Theorem:1.3} and \ref{Theorem:1.4}.

We now outline the proof of Theorem~\ref{Theorem:1.4}.
We construct the desired solution as a limit of bounded, positive classical solutions to problem~\eqref{eq:P} 
with initial data $\mu\in\mathcal{L}\cap L^\infty(\mathbb{R}^N)$ satisfying \eqref{eq:1.5} for some $T\in(0, \infty)$.
Taking $\epsilon_1\in (0,1)$ sufficiently small, 
we derive estimates for{
$$
\sup_{t\in(0,T)}|||u(t)|||_{\frac{N(p-m)}{2},\beta;\,\epsilon_1^{-1}T^{\theta}}
\quad\mbox{and}\quad 
\sup_{t\in(0,T)}\left\{t^{\frac{1}{p-1}}\|u(t)\|_{L^\infty({\mathbb R}^N)}\right\}
$$
for such a classical solution $u$ (see Proposition~\ref{Proposition:5.1}).
To this end, let $\epsilon_2\in (0,1)$ be sufficiently small, and define 
\begin{align*}
T_1:= & \,\sup\left\{t\in(0,\min\{T,T_u\})\,:\,\sup_{s\in(0,t)}|||u(s)|||_{\frac{N(p-m)}{2}, \beta;\epsilon_1^{-1}T^\theta}\le\epsilon_2\right\},\\
T_2:= & \,\sup\left\{t\in(0,\min\{T,T_u\})\,:\, \sup_{s\in(0,t)}\,
\left\{s^{\frac{1}{p-1}}\|u(s)\|_{L^\infty({\mathbb R}^N)}\right\}\le 1\right\},
\end{align*}
where $T_u$ denotes the maximal existence time of the solution~$u$.
We prove the following three estimates, by taking $\epsilon_2\in (0, 1)$ sufficiently small if necessary:
\begin{equation}
\label{eq:1.7}
\begin{split}
&\sup_{\sigma\in[\epsilon_1^{-1}t^\theta,\epsilon_1^{-1}T^\theta)}\left\{\sigma^{\frac{2}{p-m}}\sup_{z\in{\mathbb R}^N}\left(\dashint_{B(z,\sigma)}u(t)^\beta\,\dee x\right)^{\frac{1}{\beta}}\right\}
 \le C|||\mu|||_{\frac{N(p-m)}{2},\beta;\epsilon_1^{-1}T^\theta}+C\epsilon_1^{\frac{\theta'}{1-m}}	,\\
 & {|||u(t)|||_{\frac{N(p-m)}{2},\beta;\epsilon_1^{-1}t^\theta}}
 \le C\epsilon_1^{-\frac{N\theta'}{\kappa_\beta}}|||\mu|||_{\frac{N(p-m)}{2},\beta;\epsilon_1^{-1}T^\theta}^\frac{2\beta}{\kappa_\beta}
 +C\epsilon_1^{\frac{\theta'}{1-m}},\\
  &t^{\frac{1}{p-1}}\|u(t)\|_{L^\infty({\mathbb R}^N)}
 \le C\epsilon_1^{\frac{2}{p-m}-\frac{N\theta'}{\kappa_\beta}}|||\mu|||_{\frac{N(p-m)}{2},\beta;\epsilon_1^{-1}T^\theta}^\frac{2\beta}{\kappa_\beta}
 +C\epsilon_1^{\frac{2}{1-m}}, 
\end{split}
\end{equation}
for $t\in (0, \min\{T_1, T_2\})$ (see \eqref{eq:5.4} and \eqref{eq:5.6}).
We emphasize that the second terms on the right-hand side of the above estimates arise due to the assumption $m\in(0, 1)$; 
see the $L^\infty$-decay estimate (Lemma~\ref{Lemma:3.1}) and the Gronwall-type inequality for sublinear integral inequalities (Lemma~\ref{Lemma:4.1}).
Taking $\epsilon_1\in (0, 1)$ sufficiently small if necessary, we deduce from \eqref{eq:1.7} that
\begin{align*}
&|||u(t)|||_{\frac{N(p-m)}{2},\beta;\epsilon_1^{-1}T^\theta}
\le C|||\mu|||_{\frac{N(p-m)}{2},\beta;\epsilon_1^{-1}T^\theta}+C\epsilon_1^{-\frac{N\theta'}{\kappa_\beta}}|||\mu|||_{\frac{N(p-m)}{2},\beta;\epsilon_1^{-1}T^\theta}^{\frac{2\beta}{\kappa_\beta}}+\frac{\epsilon_2}{4},\\
  &t^{\frac{1}{p-1}}\|u(t)\|_{L^\infty({\mathbb R}^N)}
 \le C\epsilon_1^{\frac{2}{p-m}-\frac{N\theta'}{\kappa_\beta}}|||\mu|||_{\frac{N(p-m)}{2},\beta;\epsilon_1^{-1}T^\theta}^\frac{2\beta}{\kappa_\beta}
 +\dfrac{1}{4}, 
\end{align*}
for $t\in (0, \min\{T_1, T_2\})$ (see \eqref{eq:5.7}).
Since it follows from \eqref{eq:1.2} that 
$$
|||\mu|||_{\frac{N(p-m)}{2},\beta;\epsilon_1^{-1}T^\theta}\le C_{\epsilon_1}|||\mu|||_{\frac{N(p-m)}{2},\beta;T^\theta}
$$
for some $C_{\epsilon_1}>0$ (which is independent of $\mu$ and $T>0$),
taking $\delta>0$ sufficiently small if necessary, we deduce that 
\begin{equation}
\label{eq:1.8}
|||u(t)|||_{\frac{N(p-m)}{2},\beta;\epsilon_1^{-1}T^\theta}\le\frac{\epsilon_2}{2},
\qquad 
t^{\frac{1}{p-1}}\|u(t)\|_{L^\infty({\mathbb R}^N)}\le\frac{1}{2},
\end{equation}
for $t\in(0,\min\{T_1,T_2\})$. 
These estimates imply that $\min\{T_1, T_2\}=T$, and that $u$ satisfies \eqref{eq:1.8} for $t\in(0,T)$. 
We thus construct the desired solution as the limit of classical solutions satisfying \eqref{eq:1.8}, which enables us to complete the proof of Theorem~\ref{Theorem:1.4}.
It seems difficult to establish~\eqref{eq:1.6} without introducing the parameter $\epsilon_1$. 
The proof of Theorem~\ref{Theorem:1.3} is more complicated than that of Theorem~\ref{Theorem:1.4}, 
but it is essentially a refinement of the same strategy.

The rest of this paper is organized as follows.
In Section~2, we apply arguments similar to those used in \cite{IMS}*{Proposition~2.1} and prove Theorem~\ref{Theorem:1.1}. 
In Section~3, we derive the $L^\infty$-decay and energy estimates for solutions to problem~\eqref{eq:P}. 
Sections~4 and 5 are devoted to the proofs of Theorems~\ref{Theorem:1.3} and \ref{Theorem:1.4}, respectively. 
We also establish Corollary~\ref{Corollary:1.1}.
In Section~6, we prove Theorem~\ref{Theorem:1.2}.
%%%%%%%%%%%%%%%%%%%%%%
%%%%%%%%%%%%%%%%%%%%%%
\section{Proof of Theorem~\ref{Theorem:1.1}}
%%%%%%%%%%%%%%%%%%%%%%
%%%%%%%%%%%%%%%%%%%%%%
%
This section is devoted to the proof of Theorem~\ref{Theorem:1.1}, which relies primarily on the following proposition. 
\begin{proposition}
\label{Proposition:2.1}
Let $N\ge 1$, $m\in(0,1)$, and $p>1$. 
Let $u$ be a solution to problem~\eqref{eq:P} in ${\mathbb R}^N\times(0,T)$, where $T\in(0,\infty)$.  
Then there exists $C=C(N, m, p)>0$ such that
\begin{equation}
\label{eq:2.1}
  \left\{
  \begin{array}{ll}
  \displaystyle{\sup_{z\in{\mathbb R}^N}\mu(B(z,T^\theta))\le CT^{\theta\left(N-\frac{2}{p-m}\right)}} & \mbox{if $1<p<p_m$},\\
   \displaystyle{\sup_{z\in{\mathbb R}^N}\mu(B(z,\sigma))\le C\left[\log\left(e+\frac{T^{\theta}}{\sigma}\right)\right]^{-\frac{N}{2}}} 
    & \mbox{for $\sigma\in(0,T^\theta)$ if $p=p_m$},\vspace{5pt}\\
  \displaystyle{\sup_{z\in{\mathbb R}^N}\mu(B(z,\sigma))\le C\sigma^{N-\frac{2}{p-m}}}
    & \mbox{for $\sigma\in(0,T^\theta)$ if $p>p_m$}.
  \end{array}
\right.
\end{equation}
\end{proposition}
{\bf Proof.}
The proof is a slight modification of the arguments in \cite{IMS}*{Proposition~2.1} (see also \cite{TY}*{Section~3}). 
However, the choice of certain parameters in the proof is delicate. 
For the sake of completeness, we include the proof here 
(see also Remark~\ref{Remark:2.1}).

The proof is divided into several steps. 
Let $u$ be a solution to problem~\eqref{eq:P} in ${\mathbb R}^N\times(0,T)$, where $T\in(0,\infty)$.
\vspace{3pt}
\newline
\underline{Step 1}: 
Let $\psi\in C_c^{2;1}({\mathbb R}^N\times[0,T))$ be such that $0\le\psi\le 1$ in ${\mathbb R}^N\times[0,T)$. 
Then the following inequality holds:  
\begin{equation}
\label{eq:2.2}
\int_{\{x\in\mathbb{R}^N\,:\,\psi(x, 0)=1\}}\,\dee\mu(x)
\le C\int^T_0\int_{\mathbb{R}^N} \left(|\partial_t \psi|^{\frac{p}{p-1}}
+(|\Delta\psi|+|\nabla\psi|^2)^{\frac{p}{p-m}}\right)\,\dee x\,\dee t.
\end{equation}
The proof follows exactly the same lines as Step~1 of the proof of \cite{IMS}*{Proposition~2.1}, 
and is therefore omitted. 
\vspace{3pt}
\newline
\underline{Step 2}: 
Let $\zeta\in C^\infty(\mathbb{R})$ be such that 
$0\le \zeta\le 1$ in $\mathbb{R}$, $\zeta\equiv 1$ in $[1, \infty)$, $\zeta\equiv 0$ in $(-\infty, 0]$, 
and $|\zeta'|\le 2$ in $\mathbb{R}$. 
Let $z\in{\mathbb R}^N$, $a>0$, $\delta>0$, and set $\ell:=\lfloor 2\theta\rfloor +1$, 
where $\lfloor 2\theta\rfloor$ denotes the greatest integer less than or equal to $2\theta$.
For any $F\in C^\infty((0,\infty))$ with 
\begin{equation}
\label{eq:2.3}
F\le 0\quad\mbox{in}\quad [aT^\ell,\infty),
\end{equation} 
set 
$$
\psi(x,t):=\zeta(F(r(x,t)))\quad\mbox{with}\quad
r(x,t):=|x-z|^{\theta'\ell}+at^\ell+\delta
$$
for $(x,t)\in {\mathbb R}^N\times[0,T)$. 
Then $\psi\in C_c^{2;1}({\mathbb R}^N\times[0,\infty))$ with
$$
0\le \psi\le 1\quad\mbox{on}\quad\mathbb{R}^N\times[0, T),
\qquad 
\psi(\cdot, t)=0\quad\mbox{for}\quad t\in[T, \infty).
$$
Moreover, by the direct calculation, we obtain
\begin{alignat*}{1}
	\partial_t\psi(x, t)
	&=\zeta'(F(r))F'(r)a\ell t^{\ell-1},\\
	\nabla\psi(x, t)
	&=\zeta'(F(r))F'(r)\theta'\ell|x-z|^{\theta'\ell-2}(x-z),\\
	\Delta\psi(x, t)
	&=\zeta''(F(r))|F'(r)|^2(\theta')^2\ell^2|x-z|^{2(\theta'\ell-1)}+\zeta'(F(r))F''(r)(\theta')^2\ell^2|x-z|^{2(\theta'\ell-1)}\\
	&\quad\quad+\zeta'(F(r))F'(r)\theta'\ell(\theta'\ell+N-2)|x-z|^{\theta'\ell-2}.
\end{alignat*}
These imply that
\begin{alignat*}{1}
	|\partial_t\psi(x, t)|
	& \le Cat^{\ell-1}|F'(r)|=Ca^{\frac{1}{\ell}}(at^\ell)^{1-\frac{1}{\ell}}|F'(r)|\le Ca^\frac{1}{\ell}r^{1-\frac{1}{\ell}}|F'(r)|,
	\\ 
	|\nabla \psi(x, t)|^2
	& \le C|F'(r)|^2|x-z|^{2\theta'\ell-2}\le C|F'(r)|^2r^{2-\frac{2\theta}{\ell}},
	\\
	|\Delta \psi(x, t)| & \le C|F'(r)|^2|x-z|^{2\theta'\ell-2}+C|F''(r)||x-z|^{2\theta'\ell-2}+C|F'(r)||x-z|^{\theta'\ell-2}\\
	&\le C|F'(r)|^2r^{2-\frac{2\theta}{\ell}}+C|F''(r)|r^{2-\frac{2\theta}{\ell}}+C|F'(r)|r^{1-\frac{2\theta}{\ell}}.
\end{alignat*}
Then it follows from \eqref{eq:2.2} that
\begin{equation*}
\begin{split}
 & \sup_{z\in{\mathbb R}^N}\int_{\left\{x\in\mathbb{R}^N\,:\, F\left(|x-z|^{\theta'\ell}+\delta\right)\ge1\right\}}\,\dee \mu(x)\\
 & \le C\sup_{z\in\mathbb{R}^N}\int_{\{(x, t)\in\mathbb{R}^N\times(0, \infty)\,:\, 0\le F(r(x,t))\le 1\}} G(r(x,t))\,\dee x\,\dee t\\
 & =C\int_{\left\{(\rho, t)\in(0,\infty)\times(0, \infty)\,:\, 0\le F\left(\rho^{\theta'\ell}+at^\ell+\delta\right)\le 1\right\}} 
 \rho^{N-1}G(\rho^{\theta'\ell}+at^\ell+\delta)\,\dee\rho\,\dee t,
\end{split}
\end{equation*}
where 
$$
G(\xi):=\left(a^\frac{1}{\ell}\xi^{1-\frac{1}{\ell}}|F'(\xi)|\right)^\frac{p}{p-1}
+\left(|F'(\xi)|^2\xi^{2-\frac{2\theta}{\ell}}+|F''(\xi)|\xi^{2-\frac{2\theta}{\ell}}+|F'(\xi)|\xi^{1-\frac{2\theta}{\ell}}\right)^{\frac{p}{p-m}}.
$$
Setting $\rho=s^{\frac{2\theta}{\ell}}$ and $t=a^{-\frac{1}{\ell}}\tau^{\frac{2}{\ell}}$,
we obtain
\begin{equation}
\label{eq:2.4}
\begin{split}
 & \sup_{z\in{\mathbb R}^N}\int_{\{x\in\mathbb{R}^N\,:\, F(|x-z|^{\theta'\ell}+\delta)\ge1\}}\,\dee \mu(x)\\
 & \le
 Ca^{-\frac{1}{\ell}}\int_{\{(s, \tau)\in(0,\infty)\times(0, \infty)\,:\, 0\le F(s^2+\tau^2+\delta)\le 1\}} 
 s^{\frac{2\theta(N-1)}{\ell}}G(s^2+\tau^2+\delta)s^{\frac{2\theta}{\ell}-1}\tau^{\frac{2}{\ell}-1}\,\dee s\,\dee \tau
 \\
 &
  = Ca^{-\frac{1}{\ell}}\int_{\{\xi\in(0,\infty)\,:\, 0\le F(\xi^2+\delta)\le 1\}} \xi^{\frac{2N\theta}{\ell}+\frac{2}{\ell}-1}
 G(\xi^2+\delta)\,\dee\xi\int^\frac{\pi}{2}_0 (\cos\omega)^{\frac{2\theta}{\ell}-1} (\sin\omega)^{\frac{2}{\ell}-1}\,d\omega
 \\
 & \le Ca^{-\frac{1}{\ell}}\int_{\{\xi>0\,:\, 0\le F(\xi+\delta)\le 1\}}\xi^{\frac{N\theta}{\ell}+\frac{1}{\ell}-1} G(\xi+\delta)\,\dee \xi.
\end{split}
\end{equation}
\underline{Step 3}: 
Let $b$, $c$, and $d$ be positive constants chosen later.
Set
\[
	F(\xi):=\dfrac{1}{b}\Bigg(\log\bigg(1+\dfrac{d}{\xi}\bigg)-c\Bigg)\quad\mbox{for}\quad \xi> 0.
\]
By the direct calculation, we have
\begin{align}
\label{eq:2.5}
 & F(\xi)\ge 1\quad\mbox{if and only if}\quad \xi\le R_1:=\dfrac{d}{e^{b+c}-1},\\
\label{eq:2.6}
 & F(\xi)\ge 0\quad\mbox{if and only if}\quad \xi\le R_2:=\dfrac{d}{e^{c}-1}.
\end{align}
Furthermore, since
\begin{align*}
 & F'(\xi)=\dfrac{1}{b}\frac{1}{1+d/\xi}\bigg(-\frac{d}{\xi^2}\bigg)=-\frac{1}{b\xi}\frac{d}{\xi+d},\\
 & F''(\xi)=\dfrac{1}{b\xi^2} \frac{d}{\xi+d}+\frac{1}{b\xi} \frac{d}{(\xi+d)^2}=\frac{1}{b\xi^2}\frac{d(\xi+d)+d\xi}{(\xi+d)^2},
\end{align*}
we have 
$$
|F'(\xi)|\le \dfrac{1}{b\xi},\quad |F''(\xi)|\le \dfrac{2}{b\xi^2}.
$$
These imply that 
\begin{equation}
\label{eq:2.7}
\begin{split}
 & \Big(a^\frac{1}{\ell}\xi^{1-\frac{1}{\ell}}|F'(\xi)|\Big)^\frac{p}{p-1}
 \le a^\frac{p}{\ell(p-1)}b^{-\frac{p}{p-1}}\xi^{-\frac{p}{\ell(p-1)}},\\
 & 
 \Big(|F'(\xi)|^2\xi^{2-\frac{2\theta}{\ell}}+|F''(\xi)|\xi^{2-\frac{2\theta}{\ell}}+|F'(\xi)|\xi^{1-\frac{2\theta}{\ell}}\Big)^{\frac{p}{p-m}}
 \\
 & \qquad\quad
 \le C(b^{-\frac{2p}{p-m}}+b^{-\frac{p}{p-m}})\xi^{-\frac{2\theta p}{\ell(p-m)}}
 =C(b^{-\frac{2p}{p-m}}+b^{-\frac{p}{p-m}})\xi^{-\frac{p}{\ell(p-1)}}.
\end{split}
\end{equation}
Recalling \eqref{eq:2.3} and \eqref{eq:2.6}, we assume that 
\begin{equation}
\label{eq:2.8}
aT^\ell\ge R_2
\end{equation}
Then, applying \eqref{eq:2.4} together with \eqref{eq:2.5}, \eqref{eq:2.6}, and \eqref{eq:2.7}, 
and letting $\delta\to 0$,  
we obtain
\begin{align*}
 & \sup_{z\in{\mathbb R}^N}\mu\left(B\left(z,R_1^{\theta/\ell}\right)\right)\\
 & \le Ca^{-\frac{1}{\ell}}\big(a^\frac{p}{\ell(p-1)}b^{-\frac{p}{p-1}}+b^{-\frac{2p}{p-m}}+b^{-\frac{p}{p-m}}\big)\int^{R_2}_{R_1}\xi^{\frac{N\theta}{\ell}+\frac{1}{\ell}-1-\frac{p}{\ell(p-1)}}\,\dee \xi\\
 &=Ca^{-\frac{1}{\ell}}\big(a^\frac{p}{\ell(p-1)}b^{-\frac{p}{p-1}}+b^{-\frac{2p}{p-m}}+b^{-\frac{p}{p-m}}\big)\int^{R_2}_{R_1}\xi^{\frac{1}{\ell(p-1)}\left(\frac{N(p-m)}{2}-1\right)-1}\,\dee \xi\\
 &=Ca^{-\frac{1}{\ell}}R_1^{\frac{1}{\ell(p-1)}\left(\frac{N(p-m)}{2}-1\right)}
 \big(a^\frac{p}{\ell(p-1)}b^{-\frac{p}{p-1}}+b^{-\frac{2p}{p-m}}+b^{-\frac{p}{p-m}}\big)\int^{R_2/R_1}_{1}\xi^{\frac{1}{\ell(p-1)}\left(\frac{N(p-m)}{2}-1\right)-1}\,\dee \xi\\
 &=Ca^{-\frac{1}{\ell}}R_1^{\frac{\theta}{\ell}\left(N-\frac{2}{p-m}\right)}
 \big(a^\frac{p}{\ell(p-1)}b^{-\frac{p}{p-1}}+b^{-\frac{2p}{p-m}}+b^{-\frac{p}{p-m}}\big)
 \int^{\frac{e^{b+c}-1}{e^c-1}}_{1}\xi^{\frac{1}{\ell(p-1)}\left(\frac{N(p-m)}{2}-1\right)-1}\,\dee \xi
\end{align*}
for $a, b, c, d\in(0, \infty)$ satisfying \eqref{eq:2.8}.
Since the mapping 
$$
(0,\infty)\ni d\mapsto R:=R_1^{\frac{\theta}{\ell}}=\bigg(\dfrac{d}{e^{b+c}-1}\bigg)^{\frac{\theta}{\ell}}\in(0,\infty)
$$
is a bijection, the above estimate can be rewritten in terms of $R$ as follows:
\begin{align*}
 & \sup_{z\in{\mathbb R}^N}\mu(B(z,R))\\
 & \le Ca^{-\frac{1}{\ell}}R^{N-\frac{2}{p-m}}\big(a^\frac{p}{\ell(p-1)}b^{-\frac{p}{p-1}}+b^{-\frac{2p}{p-m}}+b^{-\frac{p}{p-m}}\big)
 \int^{\frac{e^{b+c}-1}{e^c-1}}_{1}\xi^{\frac{1}{\ell(p-1)}(\frac{N(p-m)}{2}-1)-1}\,\dee\xi
 \end{align*}
 for $a$, $b$, $c$, $R\in(0,\infty)$ with 
 $$
 aT^\ell \ge \frac{d}{e^c-1}
=\frac{e^{b+c}-1}{e^c-1}\frac{d}{e^{b+c}-1}
 =\frac{e^{b+c}-1}{e^c-1} R^{\theta'\ell},
 $$
 that is, 
 $$
 0<R\le \bigg(aT^\ell\dfrac{e^c-1}{e^{b+c}-1}\bigg)^{\frac{\theta}{\ell}}.
 $$
Letting $c\to\infty$, we obtain 
\begin{equation}
\label{eq:2.9}
\begin{split}
 & \sup_{z\in{\mathbb R}^N}\mu(B(z,R))\\
 & \le Ca^{-\frac{1}{\ell}}R^{N-\frac{2}{p-m}}\big(a^\frac{p}{\ell(p-1)}b^{-\frac{p}{p-1}}+b^{-\frac{2p}{p-m}}+b^{-\frac{p}{p-m}}\big)\int^{e^{b}}_{1}\xi^{\frac{1}{\ell(p-1)}(\frac{N(p-m)}{2}-1)-1}\,d\xi\\
\end{split}
\end{equation}
for $a$, $b$, $R\in(0,\infty)$ with $0<R\le \big(aT^\ell e^{-b}\big)^{\theta/\ell}$. 
\vspace{3pt}
\newline
\underline{Step 4}: 
We complete the proof of Proposition~\ref{Proposition:2.1}. 
Consider the case $p\neq p_m$.
Then, taking $a=b=1$ in \eqref{eq:2.9}, we have 
$$
\sup_{z\in{\mathbb R}^N}\mu(B(z,R))\le CR^{N-\frac{2}{p-m}}\quad\mbox{for}\quad 0<R\le \bigg(\dfrac{T}{e^{1/\ell}}\bigg)^\theta.
$$
This, together with \eqref{eq:1.2}, implies that
\begin{equation*}
\begin{array}{ll}
\displaystyle{\sup_{z\in\mathbb{R}^N}\mu(B(z,T^\theta))\le CT^{\theta\left(N-\frac{2}{p-m}\right)}} 
& \mbox{if}\quad p<p_m,\vspace{5pt}\\
\displaystyle{\sup_{z\in\mathbb{R}^N}\mu(B(z,\sigma))\le C\sigma^{N-\frac{2}{p-m}}} 
 & \mbox{for $\sigma\in(0,T^\theta)$ if $p>p_m$.}
\end{array}
\end{equation*}
Hence, Proposition~\ref{Proposition:2.1} follows in the case $p\neq p_m$.

It remains to consider the case $p=p_m$. 
Set $b_*:=N\ell(1-m)/2>0$. 
Then, taking $a=b^{\frac{N\ell(1-m)}{2}}$ in \eqref{eq:2.9}, 
we obtain 
\begin{equation}
\label{eq:2.10}
\sup_{z\in{\mathbb R}^N}\mu(B(z,R))
\le Ca^{-\frac{1}{\ell}}\big(a^\frac{p}{\ell(p-1)}b^{-\frac{p}{p-1}}+b^{-\frac{p}{p-m}}\big)\int^{e^{b}}_{1}\xi^{-1}\,d\xi
\le Cb^{-\frac{N}{2}}
\end{equation}
for $b\in[b_*,\infty)$ and $R\in(0,\infty)$ with 
$$
0<R\le \Big(T b^{\frac{N(1-m)}{2}}e^{-\frac{b}{\ell}}\Big)^\theta.
$$
Since the function  
$$
h(\xi):=\xi^{-\frac{N(1-m)}{2}}e^{\frac{\xi}{\ell}}
$$
is strictly increasing on $[b_*, \infty)$ and $h(\xi)\to\infty$ as $\xi\to\infty$, 
for any 
$$
R\in\left(0,AT^\theta\right)\quad\mbox{with}\quad A:=\Big(b_*^{\frac{N(1-m)}{2}}e^{-\frac{b_*}{\ell}}\Big)^\theta,
$$
we find $b_R\in[b_*,\infty)$ such that 
\begin{equation}
\label{eq:2.11}
R=\Big(T b_R^{\frac{N(1-m)}{2}}e^{-\frac{b_R}{\ell}}\Big)^\theta.
\end{equation}
Then we observe from \eqref{eq:2.10} that 
\begin{equation}
\label{eq:2.12}
\sup_{z\in{\mathbb R}^N}\mu(B(z,R))\le Cb_R^{-\frac{N}{2}}
\quad\mbox{for}\quad R\in\left(0,AT^\theta\right).
\end{equation}
 On the other hand, since $b_R\ge b_*$, by \eqref{eq:2.11}, we have
 $$
 \frac{R}{T^\theta}\ge b_*^{\frac{N\theta(1-m)}{2}}e^{-\frac{\theta}{\ell}b_R},
 $$
 that is,
\begin{equation}
\label{eq:2.13}
b_R\ge \theta'\ell\log\left(b_*^{\frac{N\theta(1-m)}{2}} \frac{T^\theta}{R}\right).
\end{equation}
 Since 
 $$
 b_*^{\frac{N\theta(1-m)}{2}} \frac{T^\theta}{R}\ge b_*^{\frac{N\theta(1-m)}{2}}A^{-1}=e^{\frac{\theta b_*}{\ell}}>1\quad\mbox{for}\quad R\in\left(0,AT^\theta\right),
 $$
 by \eqref{eq:2.12} and \eqref{eq:2.13}, we obtain 
 $$
 \sup_{z\in{\mathbb R}^N}\mu(B(z,R))\le C\left[\log\left(b_*^{\frac{N\theta(1-m)}{2}} \frac{T^\theta}{R}\right)\right]^{-\frac{N}{2}}
 \le C\left[\log\left(e+\frac{T^\theta}{R}\right)\right]^{-\frac{N}{2}}
 $$
for $R\in\left(0,AT^\theta\right)$.
This, together with \eqref{eq:1.2}, implies inequality~\eqref{eq:2.1} with $p=p_m$. 
Hence, Proposition~\ref{Proposition:2.1} follows.
$\Box$
\begin{remark}
\label{Remark:2.1}
The difference between the arguments in the proofs of Proposition~{\rm\ref{Proposition:2.1}} and \cite{IMS}*{Proposition~2.1} 
lies in whether $\ell>1$ or $\ell=1$. Indeed, in the case $m\ge 1$, we have $\theta'>2$, and hence $\ell=\lfloor2\theta\rfloor+1=1$. 
In this case, the proof of Proposition~{\rm\ref{Proposition:2.1}} is the same as that of \cite{IMS}*{Proposition~2.1} except for the argument in Step.~{\rm 4} with $p=p_m>1$.
\end{remark}
{\bf Proof of Theorem~\ref{Theorem:1.1}.}
Once Proposition~\ref{Proposition:2.1} has been established,
the proof of Theorem~\ref{Theorem:1.1} proceeds in the same way as the proof of \cite{IMS}*{Theorem~1.1}.
We omit the details.
$\Box$
%%%%%%%%%%%%%%%%%%%%%%
%%%%%%%%%%%%%%%%%%%%%%
\section{Estimates for solutions to problem~\eqref{eq:P}}
%%%%%%%%%%%%%%%%%%%%%%
%%%%%%%%%%%%%%%%%%%%%%
In this section, we derive $L^\infty$-decay and energy estimates for bounded classical solutions to equation~\eqref{eq:E}.
A function $u$ is called a bounded classical solution to equation~\eqref{eq:E} in $\mathbb{R}^N\times(0, T)$ if
$u\in C^{2; 1}(\mathbb{R}^N\times (0, T))\cap BC(\mathbb{R}^N\times (0, T))$
and satisfies \eqref{eq:E} in the classical sense. 
We begin by establishing the $L^\infty$-decay estimates.
\begin{lemma}
\label{Lemma:3.1} 
Let $N\ge 1$, $m\in(0,1)$, $p>1$, and $r\ge 1$ be such that $\kappa_r=N(m-1)+2r>0$. 
Let $u$ be a bounded classical solution to equation~\eqref{eq:E} in ${\mathbb R}^N\times(0,T)$, where $T\in(0,\infty]$. 
Then there exists $C=C(N,m,p,r)>0$ such that 
\begin{equation*}
\|u(t)\|_{L^\infty({\mathbb R}^N)}\le Ct^{-\frac{N}{\kappa_r}}\left(\sup_{s\in(0,t)}\sup_{z\in{\mathbb R}^N}\int_{B(z,R)} u(s)^r\,\dee x\right)^{\frac{2}{\kappa_r}}
+\left(\frac{t}{R^2}\right)^{\frac{1}{1-m}}
\end{equation*}
for $R>0$ and $t\in(0,T_*)$, 
where 
$$
T_*:=\sup\left\{t\in(0,T)\,:\,\sup_{s\in(0, t)}\left\{s^{\frac{1}{p-1}}\|u(s)\|_{L^\infty({\mathbb R}^N)}\right\}\le 1\right\}.
$$
\end{lemma}
{\bf Proof.} 
The proof is a slight modification of that of \cite{DK}*{Theorem~3.1}. 
Let $u$ be a bounded classical solution to equation~\eqref{eq:E} in ${\mathbb R}^N\times(0,T)$, where $T\in(0,\infty]$, 
and let $T_*$ be defined as above. 
Let $z\in{\mathbb R}^N$, $R>0$, $t\in(0,T_*)$, and $\sigma\in(0, 1)$. 
For $n=0,1,2,\dots$, set
$$
R_n:=R(1+\sigma 2^{-n}),\quad t_n:=\frac{t}{2}(1-\sigma 2^{-n}), \quad B_n:=B(z, R_n),\quad Q_n:=B_n\times(t_n,t],
$$
and let $\zeta_n\in C^\infty_{c}(Q_n)$ satisfy
\begin{align*}
 & 0\le\zeta_n\le 1\mbox{ in } Q_n,\quad
\zeta_n=1\mbox{ in }Q_{n+1},
 \quad |\nabla\zeta_n|\le\frac{2^{n+1}}{\sigma R}\mbox{ in }Q_n,\quad 0\le \partial_\tau\zeta_n\le \frac{2^{n+2}}{\sigma t}\mbox{ in }Q_n.
\end{align*}
Let $k>0$ and set 
$$
k_n:=k\left(1-\frac{1}{2^{n+1}}\right)\in\left(\frac{k}{2},k\right),
\quad
w_n:=(v-k_n)_+^{\frac{r}{2m}}\zeta_n,
$$
for $n=0, 1, 2, \dots$, where $v:=u^m$. 
Note that $v$ satisfies 
$$
\partial_t(v^{\frac{1}{m}})=\Delta v+v^{\frac{p}{m}}\quad\mbox{in}\quad Q_n.
$$
We multiply this equation by $(v-k_n)_+^{r/m-1}\zeta_n^2$, and integrate it on $B_n\times (t_n, s)$ for each $s\in(t_n, t)$.
Then
\begin{align*}
 & \int^{s}_{t_n}\int_{B_n}\partial_t(v^{\frac{1}{m}})(v-k_n)_+^{\frac{r}{m}-1}\zeta_n^2\,\dee x\,\dee \tau\\
 &\qquad=\dfrac{1}{m}\int_{B_n}\left(\int^{v(s)}_0 \xi^\frac{1-m}{m}(\xi-k_n)_+^{\frac{r}{m}-1}\,\dee \xi \right)\zeta_n(s)^2\,\dee x\\
 &\qquad\quad-\dfrac{2}{m}\int^{s}_{t_n}\int_{B_n}\left(\int^{v}_0 \xi^\frac{1-m}{m}(\xi-k_n)_+^{\frac{r}{m}-1}\,\dee \xi \right)\zeta_n\partial_\tau\zeta_n\,\dee x\,\dee \tau\\
 &\qquad\ge Ck^\frac{1-m}{m}\int_{B_n}w_n(s)^2\,\dee x-\dfrac{C2^{n}\|u\|_{L^\infty(Q_0)}^{1-m}}{\sigma t}\int^{s}_{t_n}\int_{B_n}v(v-k_n)_+^{\frac{r}{m}-1}\,\dee x\,\dee \tau,\\
 & \int^{s}_{t_n}\int_{B_n}\Delta v(v-k_n)_+^{\frac{r}{m}-1}\zeta_n^2\,\dee x\,\dee \tau\\
 &\qquad =-\dfrac{r-m}{m}\int^{s}_{t_n}\int_{B_n}|\nabla v|^2(v-k_n)_+^{\frac{r}{m}-2}\zeta_n^2\,\dee x\,\dee\tau-2\int^{s}_{t_n}\int_{B_n}\nabla v\cdot\nabla \zeta_n(v-k_n)_+^{\frac{r}{m}-1}\zeta_n\,\dee x\,\dee \tau\\
 &\qquad \le -C\int^{s}_{t_n}\int_{B_n}|\nabla w_n|^2\,\dee x\,\dee \tau+\dfrac{C2^{2n}}{\sigma^2R^2}\int^{s}_{t_n}\int_{B_n}(v-k_n)_+^{\frac{r}{m}}\,\dee x\,\dee \tau.
\end{align*}
Furthermore, it follows that
$$
v\le (v-k_{n-1})_+ +k_{n-1}
\le (v-k_{n-1})_+ +2^{n+1}(k_n-k_{n-1})
\le 2^{n+2}(v-k_{n-1})_+\quad\mbox{if}\quad v>k_n.
$$
Combining the above estimates, we obtain
\begin{align*}
&k^\frac{1-m}{m}\int_{B_n}w_n(s)^2\,\dee x+\int^{s}_{t_n}\int_{B_n}|\nabla w_n|^2\,\dee x\,\dee \tau\\
&\qquad \le \dfrac{C2^{n}\|u\|_{L^\infty(Q_0)}^{1-m}}{\sigma t}\int^{s}_{t_n}\int_{B_n}v(v-k_{n})_+^{\frac{r}{m}-1}\,\dee x\,\dee\tau+\dfrac{C2^{2n}}{\sigma^2R^2}\int^{s}_{t_n}\int_{B_n}(v-k_n)_+^{\frac{r}{m}}\,\dee x\,\dee \tau\\
&\qquad \qquad+ \int^{s}_{t_n}\int_{B_n}v^{\frac{p}{m}}(v-k_n)_+^{\frac{r}{m}-1}\,\dee x\,\dee \tau\\
&\qquad\le  \frac{C2^{2n}}{\sigma^2t}\left(\frac{t}{R^2}+\|u\|_{L^\infty(Q_0)}^{1-m}+t\|u\|_{L^\infty(Q_0)}^{p-m}\right)\int^{s}_{t_n}\int_{B_n}(v-k_{n-1})_+^{\frac{r}{m}}\,\dee x\,\dee \tau\\
&\qquad\le  \frac{C2^{2n}}{\sigma^2t}\left(\frac{t}{R^2}+\|u\|_{L^\infty(Q_0)}^{1-m}+t\|u\|_{L^\infty(Q_0)}^{p-m}\right)\int^{s}_{t_{n-1}}\int_{B_{n-1}}w_n^2\,\dee x\,\dee \tau 
\end{align*}
for $s\in(t_n, t)$, and hence,
\begin{equation}
\label{eq:3.1}
\begin{split}
 & k^{\frac{1-m}{m}}\sup_{s\in(t_n,t)}\int_{B_n}w_n(s)^2\,\dee x
+\iint_{Q_n}|\nabla w_n|^2\,\dee x\,\dee \tau\\
 & \qquad
 \le \frac{C2^{2n}}{\sigma^2t}\left(\frac{t}{R^2}+\|u\|_{L^\infty(Q_0)}^{1-m}+t\|u\|_{L^\infty(Q_0)}^{p-m}\right)\iint_{Q_{n-1}}w_{n-1}^2\,\dee x\,\dee \tau.
\end{split}
\end{equation}
If $\|u\|_{L^\infty(Q_0)}^{1-m}<t/R^2$, then the desired estimate follows immediately. 
Otherwise, we rewrite \eqref{eq:3.1} as follows:
\begin{equation*}
\begin{split}
 & k^{\frac{1-m}{m}}\sup_{s\in(t_n,t)}\int_{B_n}w_n(s)^2\,\dee x
+\iint_{Q_n}|\nabla w_n|^2\,\dee x\,\dee\tau\\
 & \qquad
 \le \frac{C2^{2n}}{\sigma^2t}\left(\|u\|_{L^\infty(Q_0)}^{1-m}+t\|u\|_{L^\infty(Q_0)}^{p-m}\right)\iint_{Q_{n-1}}w_{n-1}^2\,\dee x\,\dee\tau.
\end{split}
\end{equation*}
Then, by the definition of $T_*$, we have
\begin{equation*}
\begin{split}
 & k^{\frac{1-m}{m}}\sup_{s\in(t_n,t)}\int_{B_n}w_n(s)^2\,\dee x
+\iint_{Q_n}|\nabla w_n|^2\,\dee x\,\dee\tau\\
 & \qquad
 \le \frac{C2^{2n}}{\sigma^2t}\|u\|_{L^\infty(Q_0)}^{1-m}\iint_{Q_{n-1}}w_{n-1}^2\,\dee x\,\dee\tau.
\end{split}
\end{equation*}
Since this estimate coincides with \cite{DK}*{(3.15')}, 
applying the same argument as in the proof of \cite{DK}*{Theorem~3.1} after \cite{DK}*{(3.15')}, 
we obtain 
\begin{equation*}
\|u(t)\|_{L^\infty(B(z, R))}\le Ct^{-\frac{N}{\kappa_r}}\left(\sup_{s\in(0,t)}\int_{B(z,2R)} u(s)^r\,\dee x\,\dee s\right)^{\frac{2}{\kappa_r}}
+\left(\frac{t}{R^2}\right)^{\frac{1}{1-m}}
\end{equation*}
for $z\in{\mathbb R}^N$, $R>0$, and $t\in(0,T_*)$. 
This, together with \eqref{eq:1.2}, implies the desired estimate in Lemma~\ref{Lemma:3.1}. 
Thus, the proof is complete. 
$\Box$
\vspace{5pt}

In the following two lemmas, 
we derive energy estimates for bounded, positive classical solutions to problem~\eqref{eq:P} in the critical and supercritical cases, respectively.
A function $u$ is called a classical solution to problem~\eqref{eq:P} in $\mathbb{R}^N\times(0, T)$ with initial data $\mu \in\mathcal{L}\cap L^\infty(\mathbb{R}^N)$ 
if $u\in C^{2; 1}(\mathbb{R}^N\times(0, T))$ satisfies equation~\eqref{eq:E} in the classical sense and $u(t)\to \mu$ as $t\searrow0$ in $L^q_{loc}(\mathbb{R}^N)$ for $1\le q<\infty$.
\begin{lemma}
\label{Lemma:3.2}
Let $N\ge 1$, $m\in(0,1)$, $p=p_m>1$, $\mu\in \mathcal{L}\cap L^\infty({\mathbb R}^N)$, $T\in(0,\infty)$, and $\alpha>0$. 
Let $u$ be a positive classical solution to problem~\eqref{eq:P} in ${\mathbb R}^N\times(0,T)$ 
such that 
$$
\sup_{t\in(0,T)}\|u(t)\|_{L^\infty({\mathbb R}^N)}<\infty.
$$
Let $\Psi_\alpha$ and $\eta$ be defined as in \eqref{eq:1.1}. 
Then 
$$
\sup_{z\in{\mathbb R}^N}\int^t_0\int_{B(z,\sigma)}u^{m-1}\Psi_\alpha''(u)|\nabla u|^2\,\dee x\,\dee s<\infty
$$
for $t\in(0,T)$ and $\sigma>0$.
Furthermore, there exists $C=C(N,m,\alpha)>0$ such that
\begin{equation*}
\begin{split} 
 & \sup_{s\in(0,t]}\sup_{z\in{\mathbb R}^N}\int_{B(z,\sigma)}\Psi_\alpha(u(s))\,\dee x+\sup_{z\in{\mathbb R}^N}\int^t_0\int_{B(z,\sigma)}u^{m-1}\Psi_\alpha''(u)|\nabla u|^2\,\dee x\,\dee s\\
 & \le C\sup_{z\in{\mathbb R}^N}\int_{B(z,\sigma)}\Psi_\alpha(\mu)\,\dee x
 +C\sigma^{-2}\sup_{z\in{\mathbb R}^N}\int^t_0\int_{B(z,\sigma)}u^{m-1}\Psi_\alpha(u)\,\dee x\,\dee s\\
 & +C\sup_{z\in{\mathbb R}^N}\int_0^t\int_{B(z,\sigma)}\Psi_\alpha(u)\,\dee x\,\dee s
 +CM_\sigma[u](t)^{p-m}\sup_{z\in{\mathbb R}^N}\int^t_0\int_{B(z,\sigma)}u^{m-1}\Psi_\alpha''(u)|\nabla u|^2\,\dee x\,\dee s\
\end{split}
\end{equation*}
for $t\in(0,T)$ and $\sigma>0$ if 
$$
M_\sigma[u](t):=\sup_{s\in(0,t]}\sup_{z\in\mathbb{R}^N}\sup_{r\in(0,\sigma]}
\left\{\eta(r) \Psi_\alpha^{-1}\left(\,\dashint_{B(z,r)} \Psi_\alpha(u(s))\,\dee y\right)\right\}\le 1.
$$
\end{lemma}
\begin{lemma}
\label{Lemma:3.3}
Let $N\ge 1$, $m\in(0,1)$, $p>1$ with $p>p_m$, $\mu\in {\mathcal L}\cap L^\infty({\mathbb R}^N)$, $T\in(0,\infty]$, and $\beta>1$. 
Let $u$ be a positive classical solution to problem~\eqref{eq:P} in ${\mathbb R}^N\times(0,T)$
such that 
$$
\sup_{t\in(0,T)}\|u(t)\|_{L^\infty({\mathbb R}^N)}<\infty.
$$
Then 
$$
\sup_{z\in{\mathbb R}^N}\int^t_0\int_{B(z,\sigma)}u^{m+\beta-3}|\nabla u|^2\,\dee x\,\dee s<\infty
$$
for $t\in(0,T)$ and $\sigma>0$. 
Furthermore, there exists $C=C(N,m,p,\beta)>0$ such that 
\begin{equation*}
\begin{split} 
 & \sup_{s\in(0,t]}\sup_{z\in{\mathbb R}^N}\int_{B(z,\sigma)}u(s)^\beta\,\dee x+\sup_{z\in{\mathbb R}^N}\int^t_0\int_{B(z,\sigma)}u^{m+\beta-3}|\nabla u|^2\,\dee x\,\dee s\\
 & \le C\sup_{z\in{\mathbb R}^N}\int_{B(z,\sigma)}\mu^\beta\,\dee x\\
 & \qquad
  +C\left(1+\sup_{s\in(0,t)}|||u(s)|||^{p-m}_{\frac{N(p-m)}{2},\beta;\sigma}\right)
 \sigma^{-2}\sup_{z\in{\mathbb R}^N}\int^t_0\int_{B(z,\sigma)}u^{m+\beta-1}\,\dee x\,\dee s\\
 & \qquad
 +C\left(\sup_{s\in(0,t)}|||u(s)|||_{\frac{N(p-m)}{2},\beta;\sigma}\right)^{p-m} 
 \sup_{z\in{\mathbb R}^N}\int^t_0\int_{B(z,\sigma)}u^{m+\beta-3}|\nabla u|^2\,\dee x\,\dee s
\end{split}
\end{equation*}
for $t\in(0,T)$ and $\sigma>0$. 
\end{lemma}
The proofs are essentially the same as those of \cite{IMS}*{Lemma~3.1} and \cite{IMS}*{Lemma~3.2}, respectively, 
and are therefore omitted.
%

%
%%%%%%%%%%%%%%%%%%%%%%
%%%%%%%%%%%%%%%%%%%%%%
\section{Proof of Theorem~\ref{Theorem:1.3}}
%%%%%%%%%%%%%%%%%%%%%%
%%%%%%%%%%%%%%%%%%%%%%
In this section, we prove Theorem~\ref{Theorem:1.3} 
by establishing a sufficient condition for the existence of solutions to problem~\eqref{eq:P} in the critical case.
The key ingredient in the proof is the following proposition, 
which is derived from the $L^\infty$ decay and the energy estimate provided in Lemmas~\ref{Lemma:3.1} and \ref{Lemma:3.2}, respectively.

\begin{proposition}
\label{Proposition:4.1}
Let $N\ge 1$, $m\in(0,1)$, $p=p_m>1$, and $\alpha>0$. 
Let $\Psi_\alpha$ and $\eta$ be defined as in \eqref{eq:1.1}.
Then there exists $\delta>0$ such that
if $\mu\in \mathcal{L}\cap L^\infty({\mathbb R}^N)$ with $\mu^{-1}\in L^\infty(\mathbb{R}^N)$ satisfies
\begin{equation}
\label{eq:4.1}
|||\mu|||_{\eta,\Psi_\alpha;1}\le\delta,
\end{equation}
then problem~\eqref{eq:P} admits a unique classical solution $u$ in $\mathbb{R}^N\times(0, 1)$ satisfying
$$
\sup_{t\in(0,1)}|||u(t)|||_{\eta,\Psi_\alpha,1}
+\sup_{t\in(0,1)}\left\{t^{\frac{1}{p-1}}\left[\log\left(e+\frac{1}{t}\right)\right]^{\frac{1}{p-1}}\|u(t)\|_{L^\infty({\mathbb R}^N)}\right\}\le C_*,
$$
where $C_*$ is a constant depending only on $N$, $m$, $p$, and $\alpha$. 
\end{proposition}
For the proof, 
we begin by recalling some properties of $\Psi_\alpha$.
First, we have 
\begin{equation}
\label{eq:4.2}
\Psi_\alpha^{-1}(\xi)\asymp \xi[\log(e+\xi)]^{-\alpha}\quad\mbox{for $\xi\in[0,\infty)$}.
\end{equation}
Furthermore, for any fixed $k>0$, 
\begin{equation}
\label{eq:4.3}
\Psi_\alpha(k\xi)\asymp\Psi_\alpha(\xi),\quad
\Psi_\alpha^{-1}(k\xi)\asymp \Psi_\alpha^{-1}(\xi),\quad\mbox{for $\xi\in[0,\infty)$}.
\end{equation}
Since $\Psi_\alpha^{-1}$ is non-decreasing on $[0, \infty)$, we observe from \eqref{eq:4.3} that
\begin{equation}
\label{eq:4.4}
\Psi_\alpha^{-1}(a+b)\le \Psi_\alpha^{-1}(2\max\{a,b\})\le C(\Psi_\alpha^{-1}(a)+\Psi_\alpha^{-1}(b))
\quad\mbox{for $a$, $b\in[0,\infty)$}.
\end{equation}

Next, following the arguments in \cite{IMS}, 
we introduce a $C^1$-function $\gamma$ in $[0,1]$ by
$$
\int^{\gamma(\xi)}_0 s\eta(s)^{m-1}\,\dee s=C_\eta \xi\quad\mbox{for $\xi\in[0,1]$},
\quad\mbox{where $\displaystyle{C_\eta:=\int^1_0 s\eta(s)^{m-1}\,\dee s}$}.
$$
Then $\gamma'>0$ in $[0,1]$, $\gamma(0)=0$, and $\gamma(1)=1$. 
Since 
\begin{align*}
C_\eta\xi=\int^{\gamma(\xi)}_0 s\eta(s)^{m-1}\,\dee s & =\int^{\gamma(\xi)}_0 s^{1+N(m-1)}\left[\log\left(e+\frac{1}{s}\right)\right]^{\frac{N(m-1)}{2}}\,\dee s\\
 & \asymp \gamma(\xi)^{2+N(m-1)}\left[\log\left(e+\frac{1}{\gamma(\xi)}\right)\right]^{\frac{N(m-1)}{2}}
=\gamma(\xi)^2\eta(\gamma(\xi))^{m-1}
\end{align*}
for $\xi\in[0,1]$, we obtain
\begin{equation}
\label{eq:4.5}
\gamma(\xi)^2\eta(\gamma(\xi))^{m-1}\asymp\xi\quad\mbox{for}\quad \xi\in[0,1],
\end{equation}
and hence,
$$
\gamma(\xi)\asymp \xi^\frac{1}{\kappa}\left[\log\left(e+\dfrac{1}{\xi}\right)\right]^{\frac{N(1-m)}{2\kappa}}\quad\mbox{for}\quad \xi\in[0,1].
$$
This, together with \eqref{eq:4.5} and $p=p_m$, implies that
\begin{equation}
\label{eq:4.6}
\eta(\gamma(\xi))\asymp \xi^\frac{N}{\kappa}\left[\log\left(e+\dfrac{1}{\xi}\right)\right]^\frac{N}{\kappa}=\xi^\frac{1}{p-1}\left[\log\left(e+\dfrac{1}{\xi}\right)\right]^\frac{1}{p-1}\quad\mbox{for}\quad \xi\in[0,1].
\end{equation}

We prove a Gronwall-type inequality, which will be used in both the subcritical and supercritical cases.
\begin{lemma}
\label{Lemma:4.1}
Let $A_1, A_2, A_3\ge0$, $T>0$, and $m\in(0, 1)$.
Assume that a nonnegative function $f\in BC((0, T))$ satisfies
\begin{equation}
\label{eq:4.7}
f(t)\le A_1+A_2\int^t_0f(s)^m\,\dee s+A\int^t_0f(s)\,\dee s\quad\mbox{for}\quad t\in(0, T).
\end{equation}
Then
$$
f(t)\le e^{A_3 t}\left(A_1^{1-m}+(1-m)A_2t\right)^\frac{1}{1-m}\quad\mbox{for}\quad t\in(0, T).
$$
\end{lemma}
{\bf Proof.}
Set
$$
g(t):=A_1+A_2\int^t_0f(s)^m\,\dee s+A_3\int^t_0f(s)\,\dee s\quad\mbox{for}\quad t\in(0, T).
$$
Then it follows from \eqref{eq:4.7} that
$$
g'(t)=A_2f(t)^m+A_3f(t)\le A_2g(t)^m+A_3g(t)\quad\mbox{for}\quad t\in(0, T).
$$
This implies that 
\begin{align*}
\dfrac{\dee}{\dee t}\left(e^{-A_3t}g(t)\right)^{1-m}
 & =(1-m)\left(e^{-A_3t}g(t)\right)^{-m}\left(-A_3e^{-A_3t}g(t)+e^{-A_3t}g'(t)\right)\\
 & \le (1-m)A_2e^{-(1-m)A_3t}\le (1-m)A_2
\end{align*}
for $t\in(0, T)$. Then we obtain
$$
\left(e^{-A_3t}f(t)\right)^{1-m}\le
\left(e^{-A_3t}g(t)\right)^{1-m}\le g(0)^{1-m}+(1-m)A_2t=A_1^{1-m}+(1-m)A_2t
$$
for $t\in(0, T)$. This implies the desired inequality, and the proof is complete. 
$\Box$\vspace{5pt}

We now prove Proposition~\ref{Proposition:4.1}. 
\vspace{5pt}
\newline
\noindent
{\bf Proof of Proposition~\ref{Proposition:4.1}.}
Let $\delta\in(0,1)$ be sufficiently small, and let $\mu\in \mathcal{L}\cap L^\infty(\mathbb{R}^N)$ with $\mu^{-1}\in L^\infty(\mathbb{R}^N)$ satisfy \eqref{eq:4.1}.
By the classical theory of quasilinear parabolic equations (see e.g., \cite{LSU}*{Chapter~V}), 
there exists a unique positive classical solution $u$ to problem~\eqref{eq:P} in $\mathbb{R}^N\times(0, T_u)$ for some $T_u\in(0, \infty)$, such that
\begin{align}
\label{eq:4.8}
 & \sup_{t\in(0,T)}\|u(t)\|_{L^\infty({\mathbb R}^N)}<\infty\quad \mbox{for}\quad T\in(0, T_u),\\
\label{eq:4.9}
 & \limsup_{t\nearrow T_u}\|u(t)\|_{L^\infty({\mathbb R}^N)}=\infty.
\end{align} 
Note that problem~\eqref{eq:P} with $p=p_m$ admits no nontrivial global-in-time solutions (see Remark~\ref{Remark:1.1}~(2)). 

Let $\epsilon_1, \epsilon_2\in(0,1)$ be parameters to be fixed later, and set
\begin{align}
\label{eq:4.10}
T_1:= & \,\sup\left\{t\in(0,\min\{T_u,1\})\,:\,\sup_{s\in(0,t)}|||u(s)|||_{\eta,\Psi_\alpha;\epsilon_1^{-1}}\le\epsilon_2\right\},\\
\label{eq:4.11}
T_2:= & \,\sup\left\{t\in(0,\min\{T_u,1\})\,:\, \sup_{s\in(0,t)}\,
\left\{s^\frac{1}{p-1}\|u(s)\|_{L^\infty({\mathbb R}^N)}\right\}\le 1\right\}.
\end{align}
Throughout the proof, generic constants $C$ are independent of $\epsilon_1, \epsilon_2\in(0,1)$. 
We also use the notation $C_{\epsilon_1}$ for generic constants which may depend on $\epsilon_1$.
\vspace{3pt}
\newline
\noindent\underline{Step 1.}
We show that $\min\{T_1,T_2\}>0$.
Since it immediately follows from \eqref{eq:4.8} that $T_2>0$, it suffices to prove that $T_1>0$.
Set 
\begin{equation}
\label{eq:4.12}
c_*:=\sup_{t\in(0,T_u/2)}\|u(t)\|_{L^\infty(\mathbb{R}^N)}<\infty.
\end{equation}
Let $t_*\in(0, T_u/2)$ be sufficiently small such that $c_*^{p-1}t_*\le 1$. It follows that
\begin{equation}
\label{eq:4.13}
\sup_{s\in(0, t_*)} |||u(s)|||_{\eta, \Psi_\alpha;  t_*^{1/2\kappa}}\le \eta(t_*^\frac{1}{2\kappa})c_*. 
\end{equation}
On the other hand, by the same argument as in the proof of Lemma~\ref{Lemma:3.2} (see Step 1 of the proof of \cite{IMS}*{Lemma~3.1}), we have
\begin{align*}
 & \sup_{s\in(0, t)}\sup_{z\in\mathbb{R}^N}\int_{B(z, \sigma)}\Psi_\alpha(u(s))\,\dee x
 \le C\sup_{z\in\mathbb{R}^N}\int_{B(z, \sigma)}\Psi_\alpha(\mu)\,\dee x\\
  & \qquad
  +C\sigma^{-2}\sup_{z\in\mathbb{R}^N}\int^t_0\int_{B(z, \sigma)}u^{m-1}\Psi_\alpha(u)\,\dee x\,\dee s
  +C\sup_{z\in\mathbb{R}^N}\int^t_0\int_{B(z, \sigma)}u^{p-1}\Psi_\alpha(u)\,\dee x\,\dee s
\end{align*}
for $\sigma\in(0, \epsilon^{-1}_1)$ and $t\in(0, t_*)$.
Since
\begin{equation}
\label{eq:4.14}
\begin{split}
&\int_{B(z,\sigma)}u(s)^{m-1}\Psi_\alpha(u(s))\,\dee x\\
 & =\int_{B(z,\sigma)}u(s)^m[\log(e+u(s))]^\alpha\,\dee x\\
 & \le \left(\int_{B(z,\sigma)}[\log(e+u(s))]^\alpha\,\dee x\right)^{1-m}\left(\int_{B(z,\sigma)}\Psi_\alpha(u(s))\,\dee x\right)^m,
\end{split}
\end{equation}
we observe from \eqref{eq:4.12} that
\begin{align*}
&\sup_{s\in(0, t)}\sup_{z\in\mathbb{R}^N}\int_{B(z, \sigma)}\Psi_\alpha(u(s))\,\dee x\\
&\le C\sup_{z\in\mathbb{R}^N}\int_{B(z, \sigma)}\Psi_\alpha(\mu)\,\dee x+C\sigma^{-\kappa}\left(\log(e+c_*)\right)^{(1-m)\alpha}\int^t_0\left(\sup_{z\in\mathbb{R}^N}\int_{B(z, \sigma)}\Psi_\alpha(u)\,\dee x\right)^m\,\dee s\\
&\quad\quad +Cc^{p-1}_*\int^t_0\left(\sup_{z\in\mathbb{R}^N}\int_{B(z, \sigma)}\Psi_\alpha(u)\,\dee x\right)\,\dee s
\end{align*}
for $\sigma\in(0, \epsilon^{-1}_1)$ and $t\in(0, t_*)$.
Then we deduce from Lemma~\ref{Lemma:4.1} and $c^{p-1}_*t_*\le1$ that
\begin{align*}
&\sup_{s\in(0, t)}\sup_{z\in\mathbb{R}^N}\int_{B(z, \sigma)}\Psi_\alpha(u(s))\,\dee x\\
& \le e^{Cc_*^{p-1}t}\left(\left(C\sup_{z\in\mathbb{R}^N}\int_{B(z, \sigma)}\Psi_\alpha(\mu)\,\dee x\right)^{1-m}+(1-m)C\sigma^{-\kappa}\left(\log(e+c_*)\right)^{(1-m)\alpha}t\right)^{\frac{1}{1-m}}\\
& \le C\sup_{z\in\mathbb{R}^N}\int_{B(z, \sigma)}\Psi_\alpha(\mu)\,\dee x+C\sigma^{-\frac{\kappa}{1-m}}\left(\log(e+c_*)\right)^\alpha t^\frac{1}{1-m}
\end{align*}
for $\sigma\in(0, \epsilon^{-1}_1)$ and $t\in(0, t_*)$.
This, together with \eqref{eq:4.2}, \eqref{eq:4.3}, and \eqref{eq:4.4}, implies that
\begin{equation}
\label{eq:4.15}
\begin{split}
&\sup_{s\in(0, t)}\sup_{z\in\mathbb{R}^N}\sup_{\sigma\in(t_*^{1/2\kappa}, \epsilon_1^{-1})}
\left\{\eta(\sigma)\Psi_\alpha^{-1}\left(\dashint_{B(z, \sigma)}\Psi_\alpha(u(s))\,\dee x\right)\right\}\\
&\le C|||\mu|||_{\eta,\Psi_\alpha;\epsilon_1^{-1}}
+C\sup_{\sigma\in(t_*^{1/2\kappa}, \epsilon_1^{-1})}
\left\{\eta(\sigma)\Psi_\alpha^{-1}\left(\sigma^{-\frac{2}{1-m}}\left(\log(e+c_*)\right)^\alpha t^\frac{1}{1-m}_*\right)\right\}\\
&\le C|||\mu|||_{\eta,\Psi_\alpha;\epsilon_1^{-1}}
+C\left(\log(e+c_*)\right)^\alpha t^\frac{1}{1-m}_*\sup_{\sigma\in(t_*^{1/2\kappa}, \epsilon_1^{-1})}\left\{\eta(\sigma)\sigma^{-\frac{2}{1-m}}\right\}\\
&= C|||\mu|||_{\eta,\Psi_\alpha;\epsilon_1^{-1}}
+C\left(\log(e+c_*)\right)^\alpha t^\frac{1}{1-m}_*\sup_{\sigma\in(t_*^{1/2\kappa}, \epsilon_1^{-1})}\left\{\sigma^{-\frac{\kappa}{1-m}}\left(\log\left(e+\dfrac{1}{\sigma}\right)\right)^\frac{N}{2}\right\}\\
& \le C|||\mu|||_{\eta,\Psi_\alpha;\epsilon_1^{-1}}
+C\left(\log(e+c_*)\right)^\alpha t_*^{\frac{1}{2(1-m)}}|\log t_*|^\frac{N}{2}
\end{split}
\end{equation}
for $t\in(0, t_*)$.
Therefore, 
taking $t_*\in(0,1)$ sufficiently small if necessary, 
we deduce from \eqref{eq:4.13} and \eqref{eq:4.15} that 
\begin{equation}
\label{eq:4.16}
\sup_{s\in(0,t)}|||u(s)|||_{\eta,\Psi_\alpha;\epsilon_1^{-1}}\le C|||\mu|||_{\eta,\Psi_\alpha;\epsilon_1^{-1}}+\dfrac{\epsilon_2}{4}
\end{equation}
for $t\in(0, t_*)$.
On the other hand, 
by \eqref{eq:1.2}, \eqref{eq:4.2}, and \eqref{eq:4.3},
we have
\begin{equation*}
\begin{split}
	 & \sup_{\sigma\in[1,\epsilon_1^{-1})}
\sup_{z\in{\mathbb R}^N}\left\{\eta(\sigma)\Psi_\alpha^{-1}\left(\dashint_{B(z,\sigma)}\Psi_\alpha(\mu)\,\dee x\right)\right\}\\
 & \le C_{\epsilon_1}\sup_{z\in{\mathbb R}^N}\left\{\eta\left(\dfrac{1}{2}\right)\Psi_\alpha^{-1}\left(C_{\epsilon_1}\,\,\dashint_{B(z,1/2)}\Psi_\alpha(\mu)\,\dee x\right)\right\}
\le C_{\epsilon_1}|||\mu|||_{\eta,\Psi_\alpha;1}.
\end{split}
\end{equation*}
This implies that 
\begin{equation}
\label{eq:4.17}
\begin{split}
|||\mu|||_{\eta,\Psi_\alpha;\epsilon_1^{-1}} & \le |||\mu|||_{\eta,\Psi_\alpha;1}+\sup_{\sigma\in[1,\epsilon_1^{-1})}
\sup_{z\in{\mathbb R}^N}\left\{\eta(\sigma)\Psi_\alpha^{-1}\left(\dashint_{B(z,\sigma)}\Psi_\alpha(\mu)\,\dee x\right)\right\}\\
 & \le C_{\epsilon_1}|||\mu|||_{\eta,\Psi_\alpha;1}.
\end{split}
\end{equation}
Therefore, taking $\delta>0$ sufficiently small if necessary, we observe from \eqref{eq:4.1} and \eqref{eq:4.16} that 
$$
\sup_{s\in(0,t)}|||u(s)|||_{\eta,\Psi_\alpha;\epsilon_1^{-1}}\le \dfrac{\epsilon_2}{2}
$$
for $t\in(0, t_*)$. 
This implies that $T_2>t_*>0$, and hence, we obtain $\min\{T_1,T_2\}>0$. 
\vspace{3pt}
\newline
\noindent\underline{Step 2.}
We obtain an estimate for
$$
\sup_{z\in{\mathbb R}^N}\sup_{\sigma\in[\epsilon_1^{-1}\gamma(t),\epsilon_1^{-1})}\left\{\eta(\sigma)\Psi_\alpha^{-1}\left(\dashint_{B(z,\sigma)}\Psi_\alpha(u(t))\,\dee x\right)\right\}
$$
for $t\in(0, T_1)$.
It follows from Lemma~\ref{Lemma:3.2} and \eqref{eq:4.10} that, 
by taking $\epsilon_2>0$ sufficiently small if necessary, we have
\begin{equation}
\label{eq:4.18}
\begin{split}
 & \sup_{z\in{\mathbb R}^N}\int_{B(z,\sigma)}\Psi_\alpha(u(t))\,\dee x
 \le C\sup_{z\in{\mathbb R}^N}\int_{B(z,\sigma)}\Psi_\alpha(\mu)\,\dee x\\
 & 
 +C\sigma^{-2}\int^t_0\left(\sup_{z\in{\mathbb R}^N}\int_{B(z,\sigma)}u^{m-1}\Psi_\alpha(u)\,\dee x\right)\,\dee s
 +C\int^t_0\left(\sup_{z\in{\mathbb R}^N}\int_{B(z,\sigma)}\Psi_\alpha(u)\,\dee x\right)\,\dee s
\end{split}
\end{equation} 
for $t\in(0,T_1)$ and $\sigma\in(0,\epsilon_1^{-1})$.
Since the map $\tau\mapsto \Psi_\alpha(e^{\tau^{1/\alpha}}-e)$ is convex on $(0, \infty)$, we observe from Jensen's inequality that 
\begin{equation}
\label{eq:4.19}
\begin{split}
 \dashint_{B(z,\sigma)}[\log(e+u(s))]^\alpha\,\dee x
& \le \left[\log\left(e+\Psi_\alpha^{-1}\left(\dashint_{B(z,\sigma)}\Psi_\alpha(u(s))\,\dee x\right)\right)\right]^{\alpha}\\
 & \le \left[\log\left(e+\eta(\sigma)^{-1}|||u(s)|||_{\eta,\Psi_\alpha;\epsilon_1^{-1}}\right)\right]^{\alpha}\\
 &\le C\left[\log\left(e+\frac{1}{\sigma}\right)\right]^\alpha
\end{split}
\end{equation}
for $z\in{\mathbb R}^N$, $s\in(0,T_1)$, and $\sigma\in(0,\epsilon_1^{-1})$.
Combining \eqref{eq:4.14}, \eqref{eq:4.18}, and \eqref{eq:4.19}, we obtain
\begin{equation*}
\begin{split}
 & \sup_{z\in{\mathbb R}^N}\int_{B(z,\sigma)}\Psi_\alpha(u(t))\,\dee x
 \le C\sup_{z\in{\mathbb R}^N}\int_{B(z,\sigma)}\Psi_\alpha(\mu)\,\dee x\\
 & \qquad\quad
 +C\sigma^{-\kappa}\left[\log\left(e+\frac{1}{\sigma}\right)\right]^{(1-m)\alpha}\int^t_0\left(\sup_{z\in{\mathbb R}^N}\int_{B(z,\sigma)}\Psi_\alpha(u)\,\dee x\right)^m\,\dee s\\
 &  \qquad\quad
+C\int^t_0\left(\sup_{z\in{\mathbb R}^N}\int_{B(z,\sigma)}\Psi_\alpha(u)\,\dee x\right)\,\dee s
\end{split}
\end{equation*} 
for $t\in(0,T_1)$ and $\sigma\in(0,\epsilon_1^{-1})$.
Then, by Lemma~\ref{Lemma:4.1}, we have
\begin{align*}
&\sup_{z\in{\mathbb R}^N}\int_{B(z,\sigma)}\Psi_\alpha(u(t))\,\dee x\\
& \le e^{Ct}\left(\left(C\sup_{z\in{\mathbb R}^N}\int_{B(z,\sigma)}\Psi_\alpha(\mu)\,\dee x\right)^{1-m}+(1-m)C\sigma^{-\kappa}\left[\log\left(e+\frac{1}{\sigma}\right)\right]^{(1-m)\alpha}t\right)^\frac{1}{1-m}\\
&\le C\sup_{z\in{\mathbb R}^N}\int_{B(z,\sigma)}\Psi_\alpha(\mu)\,\dee x+C\sigma^{-\frac{\kappa}{1-m}}\left[\log\left(e+\frac{1}{\sigma}\right)\right]^{\alpha}t^\frac{1}{1-m}
\end{align*}
for $t\in(0,T_1)$ and $\sigma\in(0,\epsilon_1^{-1})$.
Since \eqref{eq:4.5} implies that
\begin{align*}
\sup_{\sigma\in [\epsilon_1^{-1}\gamma(t), \epsilon_1^{-1})}
\left\{\sigma^{-\frac{\kappa}{1-m}}\left[\log\left(e+\dfrac{1}{\sigma}\right)\right]^{\frac{N}{2}}\right\}
&=\epsilon_1^{\frac{\kappa}{1-m}}\gamma(t)^{-\frac{\kappa}{1-m}}\left[\log\left(e+\dfrac{\epsilon_1}{\gamma(t)}\right)\right]^{\frac{N}{2}}\\
&\le\epsilon_1^{\frac{\kappa}{1-m}}\gamma(t)^{-\frac{2}{1-m}}\eta(\gamma(t))\le C\epsilon_1^{\frac{\kappa}{1-m}}t^{-\frac{1}{1-m}}
\end{align*}
for $t\in(0, 1)$, we obtain
$$
\sup_{z\in{\mathbb R}^N}\int_{B(z,\sigma)}\Psi_\alpha(u(t))\,\dee x
\le C\sup_{z\in{\mathbb R}^N}\int_{B(z,\sigma)}\Psi_\alpha(\mu)\,\dee x+C\epsilon_1^{\frac{\kappa}{1-m}}\left[\log\left(e+\frac{1}{\sigma}\right)\right]^{-\frac{N}{2}+\alpha}
$$
for $t\in(0,T_1)$ and $\sigma\in[\epsilon_1^{-1}\gamma(t),\epsilon_1^{-1})$.
This, together with \eqref{eq:4.3} and \eqref{eq:4.4}, implies that 
\begin{equation}
\label{eq:4.20}
\begin{split}
&\sup_{z\in{\mathbb R}^N}\left\{\eta(\sigma)\Psi_\alpha^{-1}\left(\dashint_{B(z,\sigma)}\Psi_\alpha(u(t))\,\dee x\right)\right\}\\
&\le C|||\mu|||_{\eta,\Psi_\alpha;\epsilon_1^{-1}}+C\eta(\sigma)\Psi_\alpha^{-1}\left(\epsilon_1^{\frac{\kappa}{1-m}}\sigma^{-N}\left[\log\left(e+\frac{1}{\sigma}\right)\right]^{-\frac{N}{2}+\alpha}\right)
\end{split}
\end{equation}
for $t\in(0,T_1)$ and $\sigma\in[\epsilon_1^{-1}\gamma(t),\epsilon_1^{-1})$.
Since 
\begin{equation}
\label{eq:4.21}
\begin{split}
\log\left(e+a\right) & \le \log\left(\frac{1}{b}\left(e+ab\right)\right)
=\log\left(\frac{1}{b}\right)+\log\left(e+ab\right)\\
 & \le 2\log\left(\frac{1}{b}\right)\log\left(e+ab\right)
=2|\log b|\log\left(e+ab\right)
\end{split}
\end{equation}
for $a>0$ and $b\in(0, e^{-1})$, we deduce from \eqref{eq:4.2} that
\begin{align*}
 & \Psi_\alpha^{-1}\left(\epsilon_1^{\frac{\kappa}{1-m}}\sigma^{-N}\left[\log\left(e+\frac{1}{\sigma}\right)\right]^{-\frac{N}{2}+\alpha}\right)\\
 & \le C\epsilon_1^{\frac{\kappa}{1-m}}\sigma^{-N}\left[\log\left(e+\frac{1}{\sigma}\right)\right]^{-\frac{N}{2}+\alpha}
 \left[\log\left(e+\epsilon_1^{\frac{\kappa}{1-m}}\sigma^{-N}\left[\log\left(e+\frac{1}{\sigma}\right)\right]^{-\frac{N}{2}+\alpha}\right)\right]^{-\alpha}\\
 & \le C\epsilon_1^{\frac{\kappa}{1-m}}\sigma^{-N}\left[\log\left(e+\frac{1}{\sigma}\right)\right]^{-\frac{N}{2}+\alpha}
 |\log\epsilon_1|^\alpha\left[\log\left(e+\sigma^{-N}\left[\log\left(e+\frac{1}{\sigma}\right)\right]^{-\frac{N}{2}+\alpha}\right)\right]^{-\alpha}\\
 & \le C\epsilon_1^{\frac{\kappa}{1-m}}|\log\epsilon_1|^\alpha\sigma^{-N}\left[\log\left(e+\frac{1}{\sigma}\right)\right]^{-\frac{N}{2}}
 =C\epsilon_1^{\frac{\kappa}{1-m}}|\log\epsilon_1|^\alpha\eta(\sigma)^{-1}
\end{align*}
for $\sigma>0$. 
This implies that 
\begin{equation}
\label{eq:4.22}
\eta(\sigma)\Psi_\alpha^{-1}\left(\epsilon_1^{\frac{\kappa}{1-m}}\sigma^{-N}\left[\log\left(e+\frac{1}{\sigma}\right)\right]^{-\frac{N}{2}+\alpha}\right)
\le C\epsilon_1^{\frac{\kappa}{1-m}}|\log\epsilon_1|^\alpha 
\end{equation}
for $\sigma>0$. 
Combining \eqref{eq:4.20} and \eqref{eq:4.22}, we obtain 
\begin{equation}
\label{eq:4.23}
\begin{split}
&\sup_{z\in{\mathbb R}^N}\sup_{\sigma\in[\epsilon_1^{-1}\gamma(t),\epsilon_1^{-1})}\left\{\eta(\sigma)\Psi_\alpha^{-1}\left(\dashint_{B(z,\sigma)}\Psi_\alpha(u(t))\,\dee x\right)\right\}\\
&\le C|||\mu|||_{\eta,\Psi_\alpha;\epsilon_1^{-1}}+C\epsilon_1^{\frac{\kappa}{1-m}}|\log\epsilon_1|^\alpha,
\quad t\in(0,T_1).
\end{split}
\end{equation}
\vspace{3pt}
\newline
\noindent\underline{Step 3.}
We obtain estimates for $\|u(t)\|_{L^\infty(\mathbb{R}^N)}$ and $|||u(t)|||_{\eta, \Psi_\alpha;\epsilon_1^{-1}\gamma(t)}$ for $t\in(0, \min\{T_1, T_2\})$.
Since
$$
\sup_{\sigma\in(0,\epsilon_1^{-1}\gamma(t))}\eta(\sigma)
= \epsilon_1^{-N}\gamma(t)^N\left[\log\left(e+\frac{\epsilon_1}{\gamma(t)}\right)\right]^{\frac{N}{2}}\le\epsilon_1^{-N}\eta(\gamma(t))
$$
for $t\in(0,1)$, we have
\begin{equation}
\label{eq:4.24}
|||u(t)|||_{\eta, \Psi_\alpha;\epsilon_1^{-1}\gamma(t)}
\le \|u(t)\|_{L^\infty({\mathbb R}^N)}\sup_{\sigma\in(0,\epsilon_1^{-1}\gamma(t))}\eta(\sigma)
\le \epsilon_1^{-N}\eta(\gamma(t))\|u(t)\|_{L^\infty({\mathbb R}^N)}
\end{equation}
for $t\in(0,T_1)$. 
On the other hand, it follows from \eqref{eq:4.21} that 
\begin{align*}
\eta(\epsilon_1^{-1}\gamma(t)) & =\epsilon_1^{-N}\gamma(t)^N\left[\log\left(e+\frac{\epsilon_1}{\gamma(t)}\right)\right]^{\frac{N}{2}}\\
 & \ge C\epsilon_1^{-N}|\log\epsilon_1|^{-\frac{N}{2}}\gamma(t)^N\left[\log\left(e+\frac{1}{\gamma(t)}\right)\right]^{\frac{N}{2}}
=C\epsilon_1^{-N}|\log\epsilon_1|^{-\frac{N}{2}}\eta(\gamma(t))
\end{align*}
for $t\in(0,1)$. 
This, together with Jensen's inequality and \eqref{eq:4.23}, implies that
\begin{equation}
\label{eq:4.25}
\begin{split}
 & \sup_{s\in(0,t)}\int_{B(z,\epsilon_1^{-1}\gamma(t))}u(s)\,\dee x\\
 & \le \epsilon_1^{-N}\gamma(t)^N\eta(\epsilon_1^{-1}\gamma(t))^{-1}\sup_{s\in(0,t)}\eta(\epsilon_1^{-1}\gamma(t))\Psi_\alpha^{-1}\left(\dashint_{B(z,\epsilon_1^{-1}\gamma(t))}\Psi_\alpha(u(s))\,\dee x\right)\\
 & \le C|\log\epsilon_1|^{\frac{N}{2}}\gamma(t)^N\eta(\gamma(t))^{-1}\left(|||\mu|||_{\eta,\Psi_\alpha;\epsilon_1^{-1}}+\epsilon_1^{\frac{\kappa}{1-m}}|\log\epsilon_1|^\alpha \right)
\end{split}
\end{equation}
for $t\in(0,T_1)$.
By \eqref{eq:4.11} and \eqref{eq:4.25}, 
we apply Lemma~\ref{Lemma:3.1} to obtain 
\begin{equation*}
\begin{split}
 \|u(t)\|_{L^\infty({\mathbb R}^n)}
 & \le Ct^{-\frac{N}{\kappa}}\left(\sup_{s\in(0,t)}\int_{B(z,\epsilon_1^{-1}\gamma(t))}u(s)\,\dee x\right)^{\frac{2}{\kappa}}
+\left(\frac{t}{\epsilon_1^{-2}\gamma(t)^2}\right)^{\frac{1}{1-m}}\\
 & \le C|\log\epsilon_1|^{\frac{N}{\kappa}}t^{-\frac{N}{\kappa}}\gamma(t)^{\frac{2N}{\kappa}}\eta(\gamma(t))^{-\frac{2}{\kappa}}
 \left(|||\mu|||_{\eta,\Psi_\alpha;\epsilon^{-1}_1}^{\frac{2}{\kappa}}+\epsilon_1^{\frac{2}{1-m}}|\log\epsilon_1|^{\frac{2\alpha}{\kappa}}\right)\\
 & \qquad\qquad
 +\epsilon_1^{\frac{2}{1-m}}\left(\frac{t}{\gamma(t)^2}\right)^{\frac{1}{1-m}}\\
\end{split}
\end{equation*}
for $t\in(0,\min\{T_1,T_2\})$.
This, together with \eqref{eq:4.5}, implies that
\begin{equation}
\label{eq:4.26}
\begin{split}
&\eta(\gamma(t))\|u(t)\|_{L^\infty({\mathbb R}^N)}\\
&\le C|\log\epsilon_1|^{\frac{N}{\kappa}}t^{-\frac{N}{\kappa}}\left(\gamma(t)^2\eta(\gamma(t))^{m-1}\right)^{\frac{N}{\kappa}}\left(|||\mu|||_{\eta,\Psi_\alpha;\epsilon_1^{-1}}^{\frac{2}{\kappa}}+\epsilon_1^{\frac{2}{1-m}}|\log\epsilon_1|^{\frac{2\alpha}{\kappa}}\right)\\
&\quad+\epsilon_1^{\frac{2}{1-m}}t^\frac{1}{1-m}\left(\gamma(t)^2\eta(\gamma(t))^{m-1}\right)^{-\frac{1}{1-m}}\\
&\le C|\log\epsilon_1|^{\frac{N}{\kappa}}|||\mu|||_{\eta,\Psi_\alpha;\epsilon_1^{-1}}^{\frac{2}{\kappa}}+C\epsilon_1^{\frac{2}{1-m}}|\log\epsilon_1|^{\frac{N+2\alpha}{\kappa}}
\end{split}
\end{equation}
for $t\in(0,\min\{T_1,T_2\})$.
Furthermore, it follows from \eqref{eq:4.24} and \eqref{eq:4.26} that
\begin{equation}
\label{eq:4.27}
|||u(t)|||_{\eta, \Psi_\alpha;\epsilon_1^{-1}\gamma(t)}
\le C\epsilon_1^{-N}|\log\epsilon_1|^{\frac{N}{\kappa}}|||\mu|||_{\eta,\Psi_\alpha;\epsilon_1^{-1}}^{\frac{2}{\kappa}}+C\epsilon_1^{\frac{\kappa}{1-m}}|\log\epsilon_1|^{\frac{N+2\alpha}{\kappa}}
\end{equation}
for $t\in(0,\min\{T_1,T_2\})$. 
\vspace{3pt}
\newline
\noindent\underline{Step 4.}
We complete the proof of Proposition~\ref{Proposition:4.1}. 
Combining \eqref{eq:4.23} and \eqref{eq:4.27}, we deduce that 
\begin{equation}
\label{eq:4.28}
\begin{split}
|||u(t)|||_{\eta,\Psi_\alpha;\epsilon_1^{-1}}
 & \le C|||\mu|||_{\eta,\Psi_\alpha;\epsilon_1^{-1}}+C\epsilon_1^{-N}|\log\epsilon_1|^{\frac{2}{\kappa}}|||\mu|||_{\eta,\Psi_\alpha;\epsilon_1^{-1}}^{\frac{2}{\kappa}}\\
 & \qquad\quad
 +C\epsilon_1^{\frac{\kappa}{1-m}}|\log\epsilon_1|^\alpha+C\epsilon_1^{\frac{\kappa}{1-m}}|\log\epsilon_1|^{\frac{N+2\alpha}{\kappa}}
\end{split}
\end{equation}
for $t\in(0,\min\{T_1,T_2\})$. 
Then, taking $\epsilon_1\in(0,1)$ sufficiently small, by \eqref{eq:4.6}, \eqref{eq:4.26}, and \eqref{eq:4.28}, we obtain 
\begin{equation}
\label{eq:4.29}
\begin{split}
 |||u(t)|||_{\eta,\Psi_\alpha;\epsilon_1^{-1}}
 &\le C|||\mu|||_{\eta,\Psi_\alpha;\epsilon_1^{-1}}+C\epsilon_1^{-N}|\log\epsilon_1|^{\frac{2}{\kappa}} |||\mu|||_{\eta,\Psi_\alpha;\epsilon_1^{-1}}^{\frac{2}{\kappa}}+\frac{\epsilon_2}{4},\\
 t^\frac{1}{p-1}\|u(t)\|_{L^\infty(\mathbb{R}^N)}
 &\le t^\frac{1}{p-1}\left[\log\left(e+\dfrac{1}{t}\right)\right]^\frac{1}{p-1}\|u(t)\|_{L^\infty(\mathbb{R}^N)}\\
 &\le C\eta(\gamma(t))\|u(t)\|_{L^\infty({\mathbb R}^N)}
\le C|\log\epsilon_1|^{\frac{N}{\kappa}}|||\mu|||_{\eta,\Psi_\alpha;\epsilon_1^{-1}}^{\frac{2}{\kappa}}+\frac{1}{4},
\end{split}
\end{equation}
for $t\in(0,\min\{T_1,T_2\})$. 
Then, by \eqref{eq:4.1}, \eqref{eq:4.17}, and \eqref{eq:4.29}, taking $\delta\in(0,1)$ sufficiently small if necessary, we obtain 
\begin{equation}
\label{eq:4.30}
\begin{split}
 |||u(t)|||_{\eta,\Psi_\alpha;\epsilon_1^{-1}}
 &\le C_{\epsilon_1}|||\mu|||_{\eta,\Psi_\alpha;1}+C_{\epsilon_1}|||\mu|||_{\eta,\Psi_\alpha;1}^{\frac{2}{\kappa}}+\frac{\epsilon_2}{4}
 \le\frac{\epsilon_1}{2},\\
 t^\frac{1}{p-1}\|u(t)\|_{L^\infty({\mathbb R}^N)}
 &\le t^\frac{1}{p-1}\left[\log\left(e+\dfrac{1}{t}\right)\right]^\frac{1}{p-1}\|u(t)\|_{L^\infty(\mathbb{R}^N)}\\
 &\le C_{\epsilon_1}|\log\epsilon_1|^{\frac{N}{\kappa}}|||\mu|||_{\eta,\Psi_\alpha;1}^{\frac{2}{\kappa}}+\frac{1}{4}\le\frac{1}{2},
\end{split}
\end{equation}
for $t\in(0,\min\{T_1,T_2\})$. 
Since
$$
(0, T_u)\ni t\mapsto |||u(t)|||_{\eta,\Psi_\alpha;\epsilon_1^{-1}}\in[0, \infty)\quad\mbox{and}\quad (0, T_u)\ni t\mapsto t^\frac{1}{p-1}\|u(t)\|_{L^\infty(\mathbb{R}^N)}\in[0, \infty)
$$
are continuous, we deduce from \eqref{eq:4.30} that $\min\{T_1, T_2\}=1$ and $u$ satisfies \eqref{eq:4.30} for $t\in(0,1)$. 
Hence, Proposition~\ref{Proposition:4.1} follows.
$\Box$\vspace{5pt}

We are now in a position to complete the proof of Theorem~\ref{Theorem:1.3}. 
\vspace{5pt}
\newline
{\bf Proof of Theorem~\ref{Theorem:1.3}.} 
In view of Remark~\ref{Remark:1.1}, it suffices to consider the case $T=1$. 
Let $\delta\in(0,1)$ be sufficiently small, and assume that $\mu\in{\mathcal L}$ satisfies \eqref{eq:1.4}. 
For $n=1,2,\dots$, we set $\mu_n:=\min\{\mu, n\}+n^{-1}$.
Then $\mu_n$, $\mu_n^{-1}\in L^\infty({\mathbb R}^N)$. 
Since $\Psi_\alpha$ and $\Psi_\alpha^{-1}$ are Lipschitz continuous on $[0, \infty)$, we observe from \eqref{eq:4.3} that there exists a natural integer $n_*$ satisfying
\begin{align*}
|||\mu_n|||_{\eta, \Psi_\alpha; 1}
&\le \sup_{z\in\mathbb{R}^N}\sup_{\sigma\in(0, 1]}\left\{\eta(\sigma)\Psi_\alpha^{-1}\left(\dashint_{B(z, \sigma)}\Psi_\alpha(\mu+n^{-1})\,\dee x\right)\right\}\\
&\le C\sup_{z\in\mathbb{R}^N}\sup_{\sigma\in(0, 1]}\left\{\eta(\sigma)\Psi_\alpha^{-1}\left(\dashint_{B(z, \sigma)}\Psi_\alpha(\mu)\,\dee x+n^{-1}\right)\right\}\\
&\le C\sup_{z\in\mathbb{R}^N}\sup_{\sigma\in(0, 1]}\left\{\eta(\sigma)\Psi_\alpha^{-1}\left(\dashint_{B(z, \sigma)}\Psi_\alpha(\mu)\,\dee x\right)\right\}+Cn^{-1}\le C\delta
\end{align*}
for $n\ge n_*$.
Then, taking $\delta>0$ sufficiently small if necessary, we deduce from Proposition~\ref{Proposition:4.1} 
that there exists a classical solution $u_n$ to  problem~\eqref{eq:P} in ${\mathbb R}^N\times(0,1)$ with the initial data $\mu_n$ 
such that 
\begin{equation}
\label{eq:4.31}
\sup_{t\in(0,1)}|||u_n(t)|||_{\eta,\Psi_\alpha,1}
+\sup_{t\in(0,1)}\left\{t^{\frac{1}{p-1}}\left[\log\left(e+\frac{1}{t}\right)\right]^{\frac{1}{p-1}}\|u_n(t)\|_{L^\infty({\mathbb R}^N)}\right\}\le C_*,
\end{equation}
where $C_*$ is a constant independent of $n$. 
We then apply \cite{DGVbook}*{Theorem~B.8.1} to obtain the following result:
\begin{itemize}
  \item  
  For any compact set $K\subset{\mathbb R}^N\times(0,1)$, 
  there exist $C>0$ and $\omega\in(0,1)$ such that
  \begin{equation*}
  |u_n(x_1,t_1)-u_n(x_2,t_2)|\le C\left(|x_1-x_2|^\omega+|t_1-t_2|^{\frac{\omega}{2}}\right)
  \end{equation*}
  for $(x_1,t_1)$, $(x_2,t_2)\in K$ and $n\ge n_*$. 
\end{itemize} 
By the Arzel\`a-Ascoli Theorem and a diagonal argument, 
we extract a subsequence $\{u_{n'}\}$ of $\{u_n\}$ and obtain a locally H\"older continuous function $u$ in ${\mathbb R}^N\times (0,1)$
such that
\begin{equation}
\label{eq:4.32}
\lim_{n'\to\infty}\|u_{n'}-u\|_{L^\infty(K)}=0
\end{equation}
for any compact subset $K$ of ${\mathbb R}^N\times(0,1)$. 
Furthermore, we observe from \eqref{eq:4.31} that 
\begin{equation*}
\begin{split}
 \int_{B(z, 1)}u_n(t)^p\,\dee x
 & \le \left\|u_n(t)^{p-1}\left(\log\left(e+u_n( t)\right)\right)^{-\alpha}\right\|_{L^\infty({\mathbb R}^N)}\int_{B(z,1)}\Psi_\alpha(u_n(t))\,\dee x\\
 & \le Ct^{-1}\left[\log\left(e+\frac{1}{t}\right)\right]^{-1-\alpha}\Psi_\alpha\left(|||u_n(t)|||_{\eta,\Psi_\alpha;1}\right)\\
 & \le Ct^{-1}\left[\log\left(e+\frac{1}{t}\right)\right]^{-1-\alpha}
\end{split}	
\end{equation*}
for $z\in{\mathbb R}^N$, $t\in(0, 1)$, and $n\ge n_*$.  
Then, by Definition~\ref{Definition:1.1}-(2) and \eqref{eq:4.32}, 
we apply the Lebesgue dominated convergence theorem and conclude that
$u$ is a solution to problem~\eqref{eq:P} in ${\mathbb R}^N\times(0,1)$ with the initial data $\mu$ and that $u$ satisfies
$$
\sup_{s\in(0,1)} s^{\frac{1}{p-1}}\left[\log\left(e+\frac{1}{s}\right)\right]^{\frac{1}{p-1}}\|u(s)\|_{L^\infty({\mathbb R}^N)}+
\sup_{s\in(0,1)}|||u(s)|||_{\eta,\Psi_\alpha;1}<\infty.
$$
Hence, Theorem~\ref{Theorem:1.3} follows.
$\Box$
%%%%%%%%%%%%%%%%%%%%%%
%%%%%%%%%%%%%%%%%%%%%%
\section{Proof of Theorem~\ref{Theorem:1.4}}
%%%%%%%%%%%%%%%%%%%%%%
%%%%%%%%%%%%%%%%%%%%%%
In this section, we establish a sufficient condition for the existence of solutions to problem~\eqref{eq:P} 
in the supercritical case and prove Theorem~\ref{Theorem:1.4}. 

By an argument similar to that in the proof of Theorem~\ref{Theorem:1.3}, 
we first prove the following proposition using 
the $L^\infty$ decay and the energy estimate provided in Lemmas~\ref{Lemma:3.1} and \ref{Lemma:3.3}, respectively. 
\begin{proposition}
\label{Proposition:5.1}
Let $N\ge 1$, $m\in(0,1)$, and $p>1$ with $p>p_m$. 
Let $\beta\in(1,N(p-m)/2)$ be such that $\kappa_\beta=N(m-1)+2\beta>0$. 
Then there exist $\delta>0$ such that
if $\mu\in \mathcal{L}\cap L^\infty({\mathbb R}^N)$ with $\mu^{-1}\in L^\infty(\mathbb{R}^N)$
satisfies
\begin{equation}
\label{eq:5.1}
|||\mu|||_{\frac{N(p-m)}{2},\beta;T^\theta}\le\delta\quad\mbox{for some $T\in(0,\infty)$},
\end{equation}
then problem~\eqref{eq:P} admits a unique classical solution $u$ in $\mathbb{R}^N\times(0,T)$ satisfying
$$
\sup_{t\in(0,1)}|||u(t)|||_{\frac{N(p-m)}{2},\beta;T^\theta}
+\sup_{t\in(0,T)}\left\{t^{\frac{1}{p-1}}\|u(t)\|_{L^\infty({\mathbb R}^N)}\right\}\le C_*
$$
where $C_*$ is a constant depending only on $N$, $m$, $p$, and $\beta$. 
\end{proposition}
{\bf Proof.}
Let $\delta\in(0,1)$ be sufficiently small, and let $\mu\in \mathcal{L}\cap L^\infty(\mathbb{R}^N)$ with $\mu^{-1}\in L^\infty(\mathbb{R}^N)$ satisfy \eqref{eq:5.1}.
By the classical theory of quasilinear parabolic equations (see e.g., \cite{LSU}*{Chapter~V}), 
there exists a unique positive classical solution $u$ to problem~\eqref{eq:P} in $\mathbb{R}^N\times(0, T_u)$ for some $T_u\in(0, \infty)$ such that \eqref{eq:4.8} and \eqref{eq:4.9} hold. 
Note that, since $\mu^{-1}\in L^\infty({\mathbb R}^N)$, 
Theorem~\ref{Theorem:1.1} implies that 
problem~\eqref{eq:P} admits no nontrivial global-in-time solutions. 

Let $\epsilon_1, \epsilon_2\in(0,1)$ be parameters to be fixed later, and set
\begin{align}
\label{eq:5.2}
T_1:= & \,\sup\left\{t\in(0,\min\{T,T_u\})\,:\,\sup_{s\in(0,t)}|||u(s)|||_{\frac{N(p-m)}{2}, \beta;\epsilon_1^{-1}T^\theta}\le\epsilon_2\right\},\\
\label{eq:5.3}
T_2:= & \,\sup\left\{t\in(0,\min\{T,T_u\})\,:\, \sup_{s\in(0,t)}\,
\left\{s^{\frac{1}{p-1}}\|u(s)\|_{L^\infty({\mathbb R}^N)}\right\}\le 1\right\}.
\end{align}
By an argument similar to that in Step~1 of the proof of Proposition~\ref{Proposition:4.1}, we have $\min\{T_1,T_2\}>0$ for sufficiently small $\delta>0$. 
Throughout the proof, generic constants $C$ are independent of $\epsilon_1, \epsilon_2\in(0,1)$.
We write $C_{\epsilon_1}$ for constants that may depends on $\epsilon_1$.

We now prove that $\min\{T_1,T_2\}=1$. 
We first obtain an estimate for
$$
\sup_{z\in{\mathbb R}^N}\sup_{\sigma\in[\epsilon_1^{-1}t^\theta,\epsilon_1^{-1}T^\theta)}\left\{\sigma^{\frac{2}{p-m}}\sup_{z\in{\mathbb R}^N}\left(\dashint_{B(z,\sigma)}u(t)^\beta\,\dee x\right)^{\frac{1}{\beta}}\right\}
$$
for $t\in(0, T_1)$.
It follows from Lemma~\ref{Lemma:3.3} and \eqref{eq:5.2} that, 
by taking $\epsilon_2>0$ small enough if necessary, we have
\begin{equation*}
\begin{split}
\sup_{z\in{\mathbb R}^N}\dashint_{B(z,\sigma)}u(t)^\beta\,\dee x
 & \le C\sup_{z\in{\mathbb R}^N}\dashint_{B(z,\sigma)}\mu^\beta\dee x
 +C\sigma^{-2}\sup_{z\in{\mathbb R}^N}\int^t_0\dashint_{B(z,\sigma)}u^{m+\beta-1}\,\dee x\,\dee s\\
  & \le C\sup_{z\in{\mathbb R}^N}\dashint_{B(z,\sigma)}\mu^\beta\dee x
  +C\sigma^{-2}\sup_{z\in{\mathbb R}^N}\int^t_0\left(\dashint_{B(z,\sigma)}u^\beta\,\dee x\right)^{\frac{m+\beta-1}{\beta}}\,\dee s
\end{split}
\end{equation*} 
for $t\in(0,T_1)$ and $\sigma>0$.
Then, by Lemma~\ref{Lemma:4.1}, we obtain
\begin{align*}
\sup_{z\in{\mathbb R}^N}\dashint_{B(z,\sigma)}u(t)^\beta\,\dee x
& \le \left(\left(C\sup_{z\in{\mathbb R}^N}\dashint_{B(z,\sigma)}\mu^\beta\dee x\right)^{\frac{1-m}{\beta}}+\dfrac{(1-m)C}{\beta}\sigma^{-2}t\right)^\frac{\beta}{1-m}\\
&\le C\sup_{z\in{\mathbb R}^N}\dashint_{B(z,\sigma)}\mu^\beta\dee x+C\sigma^{-\frac{2\beta}{1-m}}t^\frac{\beta}{1-m}
\end{align*}
for $t\in(0, T_1)$ and $\sigma>0$.
This implies that 
\begin{equation}
\label{eq:5.4}
\sup_{\sigma\in[\epsilon_1^{-1}t^\theta,\epsilon_1^{-1}T^\theta)}\left\{\sigma^{\frac{2}{p-m}}\sup_{z\in{\mathbb R}^N}\left(\dashint_{B(z,\sigma)}u(t)^\beta\,\dee x\right)^{\frac{1}{\beta}}\right\}
 \le C|||\mu|||_{\frac{N(p-m)}{2},\beta;\epsilon_1^{-1}T^\theta}+C\epsilon_1^{\frac{\theta'}{1-m}}
\end{equation}
for $t\in(0,T_1)$.

We next obtain estimates for $\|u(t)\|_{L^\infty(\mathbb{R}^N)}$ and $|||u(t)|||_{N(p-m)/2,\beta; \epsilon_1^{-1}t^\theta}$. 
It follows that
\begin{equation}
\label{eq:5.5}
|||u(t)|||_{\frac{N(p-m)}{2},\beta;\epsilon_1^{-1}t^\theta}
\le\epsilon_1^{-\frac{2}{p-m}}t^{\frac{1}{p-1}}\|u(t)\|_{L^\infty({\mathbb R}^N)},\quad t\in(0,T_1).
\end{equation}
Furthermore, by \eqref{eq:5.3} and \eqref{eq:5.4}, 
applying Lemme~\ref{Lemma:3.1}, we obtain
\begin{equation*}
\begin{split}
 & \|u(t)\|_{L^\infty({\mathbb R}^N)}
\le Ct^{-\frac{N}{\kappa_\beta}}
\left(\sup_{s\in(0,t)}\int_{B(z,\epsilon_1^{-1}t^\theta)}u(s)^\beta\,\dee x\right)^{\frac{2}{\kappa_\beta}}+\left(\frac{t}{(\epsilon_1^{-1}t^\theta)^2}\right)^{\frac{1}{1-m}}\\
 & \qquad
 = Ct^{-\frac{N}{\kappa_\beta}}\left(\epsilon_1^{-1}t^\theta\right)^{\left(N-\frac{2\beta}{p-m}\right)\frac{2}{\kappa_\beta}}
\left(\left(\epsilon_1^{-1}t^\theta\right)^{\frac{2\beta}{p-m}}\sup_{s\in(0,t)}\dashint_{B(z,\epsilon_1^{-1}t^\theta)}u(s)^\beta\,\dee x\right)^{\frac{2}{\kappa_\beta}}+\epsilon_1^{\frac{2}{1-m}}t^{-\frac{1}{p-1}}\\
&  \qquad
\le Ct^{-\frac{N}{\kappa_\beta}}\left(\epsilon_1^{-1}t^\theta\right)^{\left(N-\frac{2\beta}{p-m}\right)\frac{2}{\kappa_\beta}}
\left(|||\mu|||_{\frac{N(p-m)}{2},\beta;\epsilon_1^{-1}T^\theta}+\epsilon_1^{\frac{\theta'}{1-m}}\right)^{\frac{2\beta}{\kappa_\beta}}+\epsilon_1^{\frac{2}{1-m}}t^{-\frac{1}{p-1}}\\
&  \qquad
=Ct^{-\frac{1}{p-1}}\left(\epsilon_1^{\frac{2}{p-m}-\frac{N\theta'}{\kappa_\beta}}\left(|||\mu|||_{\frac{N(p-m)}{2},\beta;\epsilon_1^{-1}T^\theta}+\epsilon_1^{\frac{\theta'}{1-m}}\right)^{\frac{2\beta}{\kappa_\beta}}+\epsilon_1^{\frac{2}{1-m}}\right)\\
&  \qquad
\le Ct^{-\frac{1}{p-1}}\left(\epsilon_1^{\frac{2}{p-m}-\frac{N\theta'}{\kappa_\beta}}|||\mu|||_{\frac{N(p-m)}{2},\beta;\epsilon_1^{-1}T^\theta}^{\frac{2\beta}{\kappa_\beta}}+\epsilon_1^{\frac{2}{1-m}}\right)
\end{split}	
\end{equation*}
for $t\in(0,\min\{T_1,T_2\})$. 
Here, we used the following relations: 
\begin{align*}
\bullet\quad
 & -\left(N-\frac{2\beta}{p-m}\right)\frac{2}{\kappa_\beta}
 =\frac{2}{(p-m)\kappa_\beta}\left(\kappa_\beta-N(p-1)\right)
 =\frac{2}{p-m}-\frac{N\theta'}{\kappa_\beta};\\
\bullet\quad
 & -\dfrac{N}{\kappa_\beta}+\left(N-\frac{2\beta}{p-m}\right)\frac{2\theta}{\kappa_\beta}
 =-\dfrac{N}{\kappa_\beta}-\left(\frac{2}{p-m}-\frac{N\theta'}{\kappa_\beta}\right)\theta=-\frac{1}{p-1};\\
\bullet\quad
 &  \dfrac{2\beta\theta'}{(1-m)\kappa_\beta}+\frac{2}{p-m}-\frac{N\theta'}{\kappa_\beta}=\frac{\theta'}{1-m}+\frac{2}{p-m}=\frac{2(p-1)+2(1-m)}{(1-m)(p-m)}=\frac{2}{1-m}.
\end{align*}
The above estimate for $\|u(t)\|_{L^\infty({\mathbb R}^N)}$, together with~\eqref{eq:5.5}, implies that
\begin{equation}
\label{eq:5.6}
\begin{split}
 & |||u(t)|||_{\frac{N(p-m)}{2},\beta;\epsilon_1^{-1}t^\theta}
 \le C\epsilon_1^{-\frac{N\theta'}{\kappa_\beta}}|||\mu|||_{\frac{N(p-m)}{2},\beta;\epsilon_1^{-1}T^\theta}^\frac{2\beta}{\kappa_\beta}
 +C\epsilon_1^{\frac{\theta'}{1-m}},\\
 & t^{\frac{1}{p-1}}\|u(t)\|_{L^\infty({\mathbb R}^N)}
 \le C\epsilon_1^{\frac{2}{p-m}-\frac{N\theta'}{\kappa_\beta}}|||\mu|||_{\frac{N(p-m)}{2},\beta;\epsilon_1^{-1}T^\theta}^\frac{2\beta}{\kappa_\beta}
 +C\epsilon_1^{\frac{2}{1-m}}, 
\end{split}
\end{equation}
for $t\in(0,\min\{T_1,T_2\})$. 

Combining \eqref{eq:5.4} and \eqref{eq:5.6}, we have 
\begin{equation}
\label{eq:5.7}
\begin{split}
 |||u(t)|||_{\frac{N(p-m)}{2},\beta;\epsilon_1^{-1}T^\theta}
  \le C|||\mu|||_{\frac{N(p-m)}{2},\beta;\epsilon_1^{-1}T^\theta}+C\epsilon_1^{-\frac{N\theta'}{\kappa_\beta}}|||\mu|||_{\frac{N(p-m)}{2},\beta;\epsilon_1^{-1}T^\theta}^{\frac{2\beta}{\kappa_\beta}}
 +C\epsilon_1^{\frac{\theta'}{1-m}}
\end{split}
\end{equation}
for $t\in(0,\min\{T_1,T_2\})$. 
Taking $\epsilon_1\in(0,1)$ sufficiently small if necessary, we deduce from \eqref{eq:5.6} and \eqref{eq:5.7} that 
\begin{equation}
\label{eq:5.8}
\begin{split}
 & |||u(t)|||_{\frac{N(p-m)}{2},\beta;\epsilon_1^{-1}T^\theta}
\le C|||\mu|||_{\frac{N(p-m)}{2},\beta;\epsilon_1^{-1}T^\theta}+C\epsilon_1^{-\frac{N\theta'}{\kappa_\beta}}|||\mu|||_{\frac{N(p-m)}{2},\beta;\epsilon_1^{-1}T^\theta}^{\frac{2\beta}{\kappa_\beta}}+\frac{\epsilon_2}{4},\\
 & t^{\frac{1}{p-1}}\|u(t)\|_{L^\infty({\mathbb R}^N)}\le C\epsilon_1^{\frac{2}{p-m}-\frac{N\theta'}{\kappa_\beta}}|||\mu|||_{\frac{N(p-m)}{2},\beta;\epsilon_1^{-1}}^\frac{2\beta}{\kappa_\beta}+\frac{1}{4},
\end{split}
\end{equation}
for $t\in(0,\min\{T_1,T_2\})$.
On the other hand, by \eqref{eq:1.2}, 
we have
\begin{equation}
\label{eq:5.9}
|||\mu|||_{\frac{N(p-m)}{2},\beta;\epsilon_1^{-1}T^\theta}\le C_{\epsilon_1}|||\mu|||_{\frac{N(p-m)}{2},\beta;T^\theta}\le C_{\epsilon_1}\delta.
\end{equation}
Therefore, by \eqref{eq:5.8} and \eqref{eq:5.9}, 
taking $\delta\in(0,1)$ sufficiently small if necessary, we obtain 
\begin{equation}
\label{eq:5.10}
\begin{split}
 & |||u(t)|||_{\frac{N(p-m)}{2},\beta;T^\theta}\le |||u(t)|||_{\frac{N(p-m)}{2},\beta;\epsilon_1^{-1}T^\theta}\le\frac{\epsilon_2}{2},\\
 & t^{\frac{1}{p-1}}\|u(t)\|_{L^\infty({\mathbb R}^N)}\le\frac{1}{2},
\end{split}
\end{equation}
for $t\in(0,\min\{T_1,T_2\})$. 
Therefore, by the same argument as in the proof of Proposition~\ref{Proposition:4.1}, we conclude that $\min\{T_1, T_2\}=1$ and
$u$ satisfies \eqref{eq:5.10} for $t\in(0,T)$. 
Hence, Proposition~\ref{Proposition:5.1} follows. 
$\Box$
\vspace{5pt}
\newline
\noindent{\bf Proof of Theorem~\ref{Theorem:1.4}.} 
Let $\delta\in(0,1)$ be sufficiently small, and assume that $\mu\in{\mathcal L}$ satisfies \eqref{eq:1.5} for some $T\in(0, \infty]$. 
Fix $M>0$ arbitrarily and set $T_M:=\min\{T, M\}\in(0, \infty)$. 
For $n=1,2,\dots$, we set $\mu_n:=\min\{\mu, n\}+n^{-1}$. 
Then $\mu_n$, $\mu_n^{-1}\in L^\infty({\mathbb R}^N)$. 
Furthermore, there exists a natural integer $n_*=n_*(M)$ such that
\begin{align*}
	|||\mu_n|||_{\frac{N(p-m)}{2}, \beta;T_M^\theta}\le |||\mu|||_{\frac{N(p-m)}{2}, \beta;T^\theta}+T_M^{\frac{1}{p-1}} n^{-1}\le 2\delta
\end{align*}
for $n\ge n_*$.
Then, taking $\delta>0$ sufficiently small if necessary, we deduce from Proposition~\ref{Proposition:5.1} that 
there exists a classical solution $u_{n,M}$ to  problem~\eqref{eq:P} in ${\mathbb R}^N\times(0,T_M)$ with the initial data $\mu_n$ 
such that 
\begin{equation}
\label{eq:5.11}
\sup_{t\in(0,T_M)}|||u_{n,M}(t)|||_{\frac{N(p-m)}{2},\beta,T_M^\theta}
+\sup_{t\in(0,T_M)}\left\{t^{\frac{1}{p-1}}\|u_{n,M}(t)\|_{L^\infty({\mathbb R}^N)}\right\}\le C_*
\end{equation}
for $n\ge n_*$, where $C_*$ is a constant independent of $T$, $M$, and $n$. 
We then apply \cite{DGVbook}*{Theorem~B.8.1} to obtain the following result:
\begin{itemize}
  \item  
  For any $M>0$ and any compact set $K\subset{\mathbb R}^N\times(0,T_M)$, 
  there exist $C>0$ and $\omega\in(0,1)$ such that
  \begin{equation*}
  |u_{n,M}(x_1,t_1)-u_{n,M}(x_2,t_2)|\le C\left(|x_1-x_2|^\omega+|t_1-t_2|^{\frac{\omega}{2}}\right)
  \end{equation*}
  for $(x_1,t_1)$, $(x_2,t_2)\in K$ and $n\ge n_*$. 
\end{itemize} 
By the Arzel\`a-Ascoli Theorem and a diagonal argument,  
we extract a subsequence $\{u_{n',M'}\}$ of $\{u_{n,M}\}$ and obtain a locally H\"older continuous function $u$ in ${\mathbb R}^N\times(0,T)$ 
such that 
$$
\lim_{n', M'\to\infty}\|u_{n',M'}-u\|_{L^\infty(K)}=0\quad\mbox{for compact set $K\subset {\mathbb R}^N\times(0,T)$}.
$$
Furthermore, we observe from \eqref{eq:5.11} that 
\begin{equation*}
\begin{split}
& \int_{B(z,\sigma)}u_{n,M}(t)^p\,\dee x
\le C\sigma^N\left(\,\dashint_{B(z,\sigma)}u_{n,M}(t)^{\beta}\,\dee x\right)^{\frac{p}{\beta}}\\
&\qquad\quad
\le C\sigma^{N-\frac{2p}{p-m}}|||u_{n,M}(t)|||_{\frac{N(p-m)}{2},\beta;T_M^\theta}^p\le C\sigma^{N-\frac{2p}{p-m}}\quad\mbox{if}\quad 1<p<\beta,\\
& \int_{B(z,\sigma)}u_{n,M}(t)^p\,\dee x
\le \|u_{n,M}(t)\|_{L^\infty(\mathbb{R}^N)}^{p-\beta}\int_{B(z,\sigma)}u_{n,M}(t)^{\beta}\,\dee x\\
&\qquad\quad
\le C t^{-\frac{p-\beta}{p-1}}\sigma^{-\frac{2\beta}{p-m}}|||u_{n,M}(t)|||_{\frac{N(p-m)}{2},\beta;T_M^\theta}^\beta
\le C t^{-\frac{p-\beta}{p-1}}\sigma^{-\frac{2\beta}{p-m}}
\quad\mbox{if}\quad p\ge\beta,
\end{split}	
\end{equation*}
for $z\in{\mathbb R}^N$, $M>0$, $\sigma\in(0,T_M^\theta)$, $t\in(0, T_M)$, and $n\ge n_*$. 
Then, by Definition~\ref{Definition:1.1}-(2) and \eqref{eq:5.11},
we apply the Lebesgue dominated convergence theorem 
and conclude that 
$u$ is a solution to problem~\eqref{eq:P} in ${\mathbb R}^N\times(0,T)$ with initial data $\mu$ and that $u$ satisfies
$$
\sup_{t\in(0,T)}|||u(t)|||_{\frac{N(p-m)}{2},\beta;T^{\theta}}
+\sup_{t\in(0,T)}t^{\frac{1}{p-1}}\|u(t)\|_{L^\infty({\mathbb R}^N)}<\infty. 
$$
Thus, $u$ is the desired solution to problem~\eqref{eq:P}, and hence, Theorem~\ref{Theorem:1.4} follows.
$\Box$
\vspace{5pt}
\newline
{\bf Proof of Corollary~\ref{Corollary:1.1}.}
Applying the same arguments as in the proof of \cite{IMS}*{Corollary~1.1}, 
by Theorems~\ref{Theorem:1.3} and \ref{Theorem:1.4},
we obtain the desired conclusion. 
Hence, Corollary~\ref{Corollary:1.1} follows.
$\Box$
%%%%%%%%%%%%%%%%%%%%%%
%%%%%%%%%%%%%%%%%%%%%%
\section{Proof of Theorem~\ref{Theorem:1.2}}
%%%%%%%%%%%%%%%%%%%%%%
%%%%%%%%%%%%%%%%%%%%%%
In this section, we modify the arguments in Sections~4 and 5 to prove Theorem~\ref{Theorem:1.2}.
\vspace{5pt}
\newline
{\bf Proof of Theorem~\ref{Theorem:1.2}.}
By Remark~\ref{Remark:1.1}-(1), 
it suffices to consider the case of $T=1$. 
Let $1<p<p_m$, $\mu\in{\mathcal M}$, and $\delta>0$ be sufficiently small. 
Assume that condition~\eqref{eq:1.3} with $T=1$ holds.
For any $n=1, 2, \dots$, set
$$
\mu_{n}(x):=\int_{{\mathbb R}^N}\rho_{n}(x-y)\,\dee\mu (y)+n^{-1}\quad \mbox{for $x\in{\mathbb R}^N$},
$$
where $\rho_{n}$ is a mollifier in ${\mathbb R}^N$, that is, $\rho_{n}(x):=n^{N}\rho(nx)$ with $\rho\in C^\infty_{c}(B(0, 1))$ satisfying $0\le \rho\le 1$ and $\|\rho\|_{L^1(B(0, 1))}=1$.
By \eqref{eq:1.2} and \eqref{eq:1.3}, we see that $\mu_n\in BUC({\mathbb R}^N)$ for $n=1,2,\dots$, and 
we find an integer $n_*$ such that
\begin{equation*}
\begin{split}
\sup_{z\in{\mathbb R}^N}\mu_{n}(B(z, 1))
&=\sup_{z\in{\mathbb R}^N}\int_{{\mathbb R}^N}\int_{B(z, 1)}\rho_{n}(x-y)\,\dee x\,\dee\mu (y)+Cn^{-1}\\
&=\sup_{z\in{\mathbb R}^N}\int_{\mathbb{R}^N}\int_{B(z-y, 1)}\rho_{n}(w)\,\dee w\,\dee\mu (y)+Cn^{-1}\\
&=\sup_{z\in{\mathbb R}^N}\int_{\mathbb{R}^N}\rho_{n}(w)\mu(B(z-w,1))\,\dee w+Cn^{-1}\\
&\le \sup_{z\in{\mathbb R}^N}\mu(B(z, 1))+Cn^{-1}\le 2\delta
\end{split}
\end{equation*}
for $n\ge n_*$.
Furthermore, for any $n=1,2,\dots$, we can find a unique positive classical solution $u_{n}$ to problem~\eqref{eq:P} in ${\mathbb R}^N\times(0,T_{n})$ 
for some $T_n\in(0,\infty)$
such that 
\begin{align*}
 & \sup_{t\in(0, T)}\|u_{n}(t)\|_{L^\infty({\mathbb R}^N)}<\infty\quad\mbox{for $T\in(0,T_{n})$},\\
 & \limsup_{t\nearrow T_{n}}\|u_{n}(t)\|_{L^\infty({\mathbb R}^N)}=\infty,
\end{align*}
where $T_{n}$ is the maximal existence time of $u_{n}$. 
Note that problem~\eqref{eq:P} with $1<p<p_m$ admits no nontrivial global-in-time solutions (see Remark~\ref{Remark:1.1}~(2)).

Let $\epsilon_1, \epsilon_2 \in (0,1)$ be parameters to be fixed later. 
For any $n\ge n_*$, we set
$$
X_n(t):=\sup_{s\in(0, t)}\sup_{z\in\mathbb {R}^N}\int_{B(z, \epsilon_1^{-1})}u_n(s)\,\dee x,
$$
and define 
\begin{align}
\label{eq:6.1}
T^1_n:= & \,\sup\left\{t\in(0,\min\{T_n,1\})\,:\,\sup_{s\in(0,t)}X_n(s)\le\epsilon_2\right\},\\
\label{eq:6.2}
T^2_n:= & \,\sup\left\{t\in(0,\min\{T_n,1\})\,:\, \sup_{s\in(0,t)}\,
\left\{s^\frac{1}{p-1}\|u_n(s)\|_{L^\infty({\mathbb R}^N)}\right\}\le 1\right\}.
\end{align}
By an argument similar to that in Step.1 of the proof of Proposition~\ref{Proposition:4.1}, we have $\min\{T^1_n,T^2_n\}>0$ for sufficiently small $\delta>0$. 
Throughout the proof, generic constants $C$ are independent of $\epsilon_1, \epsilon_2\in(0,1)$. 
We also use the notation $C_{\epsilon_1}$ for generic constants which may depends on $\epsilon_1$.

We now prove that $\min\{T^1_n, T^2_n\}=1$ for $n \ge n_*$, provided $\delta>0$ is chosen sufficiently small.
By Lemma~\ref{Lemma:3.1} and \eqref{eq:6.2}, we see that
\begin{equation}
\label{eq:6.3}
\|u_n(t)\|_{L^\infty(\mathbb{R}^N)}\le C t^{-\frac{N}{\kappa}}X_n(t)^\frac{2}{\kappa}+\epsilon_1^\frac{2}{1-m}t^\frac{1}{1-m}
\end{equation}
for $n\ge n_*$ and $t\in (0, \min\{T^1_n, T^2_n\})$.
Let $z\in{\mathbb R}^N$ and $\zeta\in C^\infty_{c}({\mathbb R}^N)$ be such that
$$
\zeta\equiv1\ \mbox{in}\ B(z,\epsilon^{-1}_1),\quad \zeta\equiv 0\ \mbox{in}\ \mathbb{R}^N\setminus B(z,2\epsilon^{-1}_1),\quad 0\le \zeta\le 1\ \mbox{in}\ \mathbb{R}^N,\quad  
\|\nabla^2\zeta\|_{L^\infty(\mathbb{R}^N)}\le C\epsilon_1^2.
$$
Let $t\in(0, \min\{T^1_n, T^2_n\})$.
We multiply the equation $\partial_t u_n-\Delta u_n^m=u_n^p$ by $\zeta$ and integrate it on $\mathbb{R}^N\times (0, t)$.
Then we observe from \eqref{eq:1.2}, \eqref{eq:6.1}, \eqref{eq:6.2}, and \eqref{eq:6.3} that
 \begin{align*}
&\sup_{z\in{\mathbb R}^N}\int_{B(z, \epsilon_1^{-1})}u_n(t)\,\dee x\\
&\le C\sup_{z\in\mathbb{R}^N}\int_{B(z, \epsilon_1^{-1})}\mu_n\,\dee x
+C\sup_{z\in\mathbb{R}^N}\int^t_0\int_{B(z, \epsilon_1^{-1})}(\epsilon_1^{2}u_n^m+u_n^p)\,\dee x\,\dee s\\
&\le C\sup_{z\in\mathbb{R}^N}\int_{B(z, \epsilon_1^{-1})}\mu_n\,\dee x
+C\epsilon_1^{\kappa}\int^t_0X_n(s)^{m}\,\dee s\\
&\quad\quad+C\int^t_0s^{-\frac{N(p-1)}{\kappa}}X_n(s)^{\frac{2(p-1)}{\kappa}+1}\,\dee s
+C\epsilon^{\frac{2(p-1)}{1-m}}_1\int^t_0s^\frac{p-1}{1-m}X_n(s)\,\dee s\\
& \le C\sup_{z\in\mathbb{R}^N}\int_{B(z, \epsilon_1^{-1})}\mu_n\,\dee x
+C\epsilon_1^{\kappa}\int^t_0X_n(s)^{m}\,\dee s+C\epsilon_2^\frac{2(p-1)}{\kappa} X_n(t)
+C\epsilon^{\frac{2(p-1)}{1-m}}_1X_n(t)
\end{align*}
for $n\ge n_*$. Here, we used the relation that
$$
\frac{N(p-1)}{\kappa}<\frac{N(p_m-1)}{\kappa}=1,
$$
which follows from the assumption~$1<p<p_m$. 
Taking $\epsilon_1$, $\epsilon_2\in(0,1)$ sufficiently small if necessary, we have
$$
X_n(t)\le C\sup_{z\in\mathbb{R}^N}\int_{B(z, \epsilon_1^{-1})}\mu_n\,\dee x  +C\epsilon_1^{\kappa}\int^t_0X_n(s)^{m}\,\dee s
$$
for $n\ge n_*$ and $t\in (0,\min\{T^1_n, T^2_n\})$.
This, together with Lemma~\ref{Lemma:4.1}, implies that
\begin{equation*}
\begin{split}
X_n(t)
 & \le \left(\left(C\sup_{z\in\mathbb{R}^N}\int_{B(z, \epsilon_1^{-1})}\mu_n\,\dee x\right)^{1-m}+(1-m)C\epsilon_1^{\kappa}t\right)^{\frac{1}{1-m}}\\
 & \le C\sup_{z\in\mathbb{R}^N}\int_{B(z, \epsilon_1^{-1})}\mu_n\,\dee x +C\epsilon_1^\frac{\kappa}{1-m}
\end{split}
\end{equation*}
for $n\ge n_*$ and $t\in (0,\min\{T^1_n, T^2_n\})$.
Then, taking $\epsilon_1\in(0,1)$ sufficiently small if necessary, 
we obtain 
\begin{equation}
\label{eq:6.4}
X_n(t)\le C\sup_{z\in\mathbb{R}^N}\int_{B(z, \epsilon_1^{-1})}\mu_n\,\dee x +\frac{\epsilon_2}{4}
\end{equation}
for $n\ge n_*$ and $t\in (0, \min\{T^1_n, T^2_n\})$. 
Furthermore, taking $\epsilon_1$, $\epsilon_2\in(0,1)$ sufficiently small again if necessary, 
we observe from \eqref{eq:6.3} and \eqref{eq:6.4} that
\begin{equation}
\label{eq:6.5}
\begin{split}
t^\frac{1}{p-1}\|u_n(t)\|_{L^\infty(\mathbb{R}^N)}
&\le C\left\{t^{\frac{1}{p-1}\left(1-\frac{N(p-m)}{2}\right)}X_n(t)\right\}^\frac{2}{\kappa}+C\epsilon_1^\frac{2}{1-m}\\
&\le C\left(\sup_{z\in\mathbb{R}^N}\int_{B(z, \epsilon_1^{-1})}\mu_n\,\dee x\right)^{\frac{2}{\kappa}}+C\epsilon_2^\frac{2}{\kappa}+C\epsilon_1^\frac{2}{1-m}\\
& \le C\left(\sup_{z\in\mathbb{R}^N}\int_{B(z, \epsilon_1^{-1})}\mu_n\,\dee x\right)^{\frac{2}{\kappa}}+\frac{1}{4}
\end{split}
\end{equation}
for $n\ge n_*$ and $t\in (0, \min\{T^1_n, T^2_n\})$.
On the other hand, by \eqref{eq:1.2}, we have
\begin{equation}
\label{eq:6.6}
\sup_{z\in\mathbb{R}^N}\int_{B(z, \epsilon_1^{-1})}\mu_n\,\dee x \le C_{\epsilon_1}\sup_{z\in\mathbb{R}^N}\int_{B(z, 1)}\mu_n\,\dee x 
\end{equation}
for $n\ge n_*$.
Then, taking $\delta\in(0,1)$ sufficiently small if necessary, 
we deduce from \eqref{eq:6.4}, \eqref{eq:6.5}, and \eqref{eq:6.6} that 
$$
X_n(t)\le\frac{\epsilon_2}{2}\quad\mbox{and}\quad
t^\frac{1}{p-1}\|u_n(t)\|_{L^\infty(\mathbb{R}^N)}\le\frac{1}{2}
$$
for $n\ge n_*$ and $t\in (0,\min\{T^1_n, T^2_n\})$.
This, together with the same argument as in the proof of Proposition~\ref{Proposition:4.1}, we conclude from \eqref{eq:6.1} and \eqref{eq:6.2} that $\min\{T^1_n, T^2_n\}=1$ for $n\ge n_*$.
Moreover, it follows from \eqref{eq:6.3} that 
\begin{equation}
\label{eq:6.7}
X_n(t)\le\frac{\epsilon_2}{2}\quad\mbox{and}\quad
t^\frac{N}{\kappa}\|u_n(t)\|_{L^\infty(\mathbb{R}^N)}\le\frac{1}{2},
\end{equation}
for $n\ge n_*$ and $t\in (0,1)$.
We then apply \cite{DGVbook}*{Theorem~B.8.1} to obtain the following result:
\begin{itemize}
  \item  
  For any compact set $K\subset{\mathbb R}^N\times(0,1)$, 
  there exist $C>0$ and $\omega\in(0,1)$ such that
  \begin{equation*}
  |u_n(x_1,t_1)-u_n(x_2,t_2)|\le C\left(|x_1-x_2|^\omega+|t_1-t_2|^{\frac{\omega}{2}}\right)
  \end{equation*}
  for $(x_1,t_1)$, $(x_2,t_2)\in K$ and $n\ge n_*$. 
\end{itemize} 
By Arzel\`a-Ascoli's Theorem and a diagonal argument,  
we extract a subsequence $\{u_{n'}\}$ of $\{u_n\}$ and obtain a locally H\"older continuous function $u$ in ${\mathbb R}^N\times(0,1)$ 
such that 
$$
\lim_{n'\to\infty}\|u_{n'}-u\|_{L^\infty(K)}=0
$$
for any compact subset $K$ of ${\mathbb R}^N\times(0,1)$. 
Furthermore, we observe from \eqref{eq:6.7} that 
\begin{equation*}
\begin{split}
 \int_{B(z, 1)}u_{n}(t)^p\,\dee x
 & \le \|u_n(t)\|_{L^\infty(\mathbb{R}^N)}^{p-1}\int_{B(z,1)}u_n(t)\,\dee x\\
 & \le Ct^{-\frac{N(p-1)}{\kappa}}=Ct^{\frac{2}{\kappa}(1-\frac{N(p-m)}{2})-1}
\end{split}	
\end{equation*}
for $z\in{\mathbb R}^N$, $t\in(0, 1)$, and $n\ge n_*$.
Then, by Definition~\ref{Definition:1.1}-(2) and \eqref{eq:6.7}, 
we apply Lebesgue's dominated convergence theorem and conclude that
$u$ is a solution to problem~\eqref{eq:P} with the initial data $\mu$ in ${\mathbb R}^N\times(0,1)$ 
and that $u$ satisfies
$$
\sup_{t\in(0,1)} t^{\frac{N}{\kappa}}\|u(t)\|_{L^\infty({\mathbb R}^N)}+
\sup_{t\in(0,1)}\sup_{z\in\mathbb{R}^N}\int_{B(z, 1)}u(t)\,\dee x<\infty.
$$
Thus, assertion~(1) follows. 
Then assertion~(2) immediately follows from Theorem~\ref{Theorem:1.1},  \eqref{eq:1.2}, and assertion~(1). 
Hence, Theorem~\ref{Theorem:1.2} follows.
$\Box$
\vspace{8pt}

\noindent
{\bf Acknowledgment.}
K. I. was supported in part by JSPS KAKENHI Grant Number 25H00591.  
N. M. was supported in part by JSPS KAKENHI Grant Number 24K16944.

%%%%%%%%%%%%%%%%%%%%%%%%%%%%%%%%%%%%%
%%%%%%%%%%%%    references    %%%%%%%%%%%%%%%%%
%%%%%%%%%%%%%%%%%%%%%%%%%%%%%%%%%%%%%

%%%%%%%%%%%%%%%%%%%%%%
%%%%%%%%%%%%%%%%%%%%%%
\end{document}